\documentstyle[11pt]{article}
\normalbaselines
\begin{document}
\bibliographystyle{unsrt}

\def\bea*{\begin{eqnarray*}}
\def\eea*{\end{eqnarray*}}
\def\ba{\begin{array}}
\def\ea{\end{array}}
\count1=1
\def\be{\ifnum \count1=0 $$ \else \begin{equation}\fi}
\def\ee{\ifnum\count1=0 $$ \else \end{equation}\fi}
\def\ele(#1){\ifnum\count1=0 \eqno({\bf #1}) $$ \else \label{#1}\end{equation}\fi}
\def\req(#1){\ifnum\count1=0 {\bf #1}\else \ref{#1}\fi}
\def\bea(#1){\ifnum \count1=0   $$ \begin{array}{#1}
\else \begin{equation} \begin{array}{#1} \fi}
\def\eea{\ifnum \count1=0 \end{array} $$
\else  \end{array}\end{equation}\fi}
\def\elea(#1){\ifnum \count1=0 \end{array}\label{#1}\eqno({\bf #1}) $$
\else\end{array}\label{#1}\end{equation}\fi}
\def\cit(#1){
\ifnum\count1=0 {\bf #1} \cite{#1} \else 
\cite{#1}\fi}
\def\bibit(#1){\ifnum\count1=0 \bibitem{#1} [#1    ] \else \bibitem{#1}\fi}
\def\ds{\displaystyle}
\def\hb{\hfill\break}
\def\comment#1{\hb {***** {\em #1} *****}\hb }

\newcommand{\TZ}{\hbox{\bf T}}
\newcommand{\MZ}{\hbox{\bf M}}
\newcommand{\ZZ}{\hbox{\bf Z}}
\newcommand{\NZ}{\hbox{\bf N}}
\newcommand{\RZ}{\hbox{\bf R}}
\newcommand{\CZ}{\,\hbox{\bf C}}
\newcommand{\PZ}{\hbox{\bf P}}
\newcommand{\QZ}{\hbox{\bf Q}}
\newcommand{\HZ}{\hbox{\bf H}}
\newcommand{\EZ}{\hbox{\bf E}}
\newcommand{\GZ}{\,\hbox{\bf G}}
\newcommand{\DZ}{\, \hbox{\bf D}}

\newtheorem{theorem}{Theorem}
\newtheorem{lemma}{Lemma}
\newtheorem{proposition}{Proposition}
\newtheorem{corollary}{Corollary}

\vbox{\vspace{38mm}}
\begin{center}
{\LARGE \bf  Rational Curves in Rigid Calabi-Yau Three-folds \footnote{Talk presented in BIRS Workshop: Calabi-Yau Varieties and Mirror Symmetry  (December 6-11, 2003) at Banff Centre, Alberta, Canada.}}\\[10 mm] 
Shi-shyr Roan
\\{\it Institute of Mathematics, Academia Sinica \\ 
Taipei , Taiwan \\ (e-mail: maroan@ccvax.sinica.edu.tw)} 
\\[35mm]
\end{center}

\begin{abstract} We determine all the Kummer-surface-type Calabi-Yau (CY) 3-folds, i.e., those $\widehat{T/G}$ which are resolutions of 3-torus-orbifolds $T/G$ with only isolated singularities. There are only two such CY spaces: one with $G= \ZZ_3$ and $T$ being the triple-product of 1-torus carrying an order 3 automorphism, the other with $G= \ZZ_7$ and $T$ being the Jacobian of Klein quartic curve. These CY 3-folds $\widehat{T/G}$ are all rigid, hence no complex structure deformation for each of these two varieties. We further investigate problems of $\PZ^1$-curves $C$ in $\widehat{T/G}$ not contained in exceptional divisors, by considering the counting number $d$ of elements in $C$ meeting exceptional divisors in a certain manner. We have obtained the constraint of $d$. With the smallest number $d$, the complete solution of $C$ in $\widehat{T/G}$ is obtained for both cases. In the case $G=\ZZ_3$, we have derived an effective method of constructing $C$ in $\widehat{T/G}$, and obtained the explicit forms of rational curves for some other $d$ by this procedure.  

\end{abstract}

\par \vspace{5mm} \noindent

\vfill
\eject

\section{Introduction}
For the past decade, the problem of rational curves in Calabi-Yau (CY) 3-folds has drawn considerable attention of algebraic geometers, due to the bold conjecture suggested by string theorists on counting rational curves in a certain type CY 3-fold by some method using the variation of Hodge structures over complex moduli space of the corresponding mirror CY 3-fold  (see, e.g., \cite{Y}). Even though the problem has remained one theoretical challenge by now, a satisfactory solution by the rigorous mathematical reasoning should still be demanded for the answer due to explanation of the mathematical essence of CY-mirror phenomenon. Nevertheless, contrast to the role of Kummer surface in K3 surfaces, the analogous construction of CY 3-folds through torus-orbifolds often leads to rigid varieties, i.e., there is no deformation of the complex structure. Hence the counting of rational curves through the mirror-mechanism fails on those rigid CY 3-folds. Therefore, rational curves on rigid CY 3-folds would deserve certain attention for a special treatment.  
It is the objective of this paper to study problems on rational curves in the rigid CY 3-folds of Kummer-surface type, which involve abelian varieties of CM-type in number theory.

In this work, we first derive the complete structure of all CY 3-folds obtained by resolutions of 3-torus $T$ quotiented by a finite (Lie)-automorphism group $G$ with $T/G$ possessing only isolated singularities, (a study of such kind was made before in \cite{RY} under a stronger assumption (with $G$ abelian) with a weaker conclusion (of modular the isogeny-equivalence) obtained). Up to isomorphisms, there are only two such CY 3-folds $\widehat{T/G}$, and both are rigid (see Theorem \ref{thm:RCY} of the paper). One arises from $T$ being the triple-product of 1-torus carrying an order 3 automorphism with $G= \ZZ_3$, and the other is obtained by the Jacobian variety $T$ of Klein quartic curve with $G= \ZZ_7$. The crepant resolution of $T/G$ is obtained by standard methods in toric geometry  with a simple rational-surface structure of exceptional divisors for both cases. Note that these two 3-folds were also characterized as the CY 3-folds of type $III_0$ in \cite{O}, where the Lie-assumption of the group $G$ can be replaced by the biregular ones with a suitable constraint on elements in $G$ and the topology of $T/G$.   
As the exceptional divisors in these two $\widehat{T/G}$ all consist of one-parameter family of rational curves, the interesting ones among rational curves in $\widehat{T/G}$ should be those not contained in any exceptional divisor, which will be one main topic we discuss in this paper. We propose a modest and seemingly innocuous constraint on rational curves (the condition $(P)$ in Section 3), which, though not cover all, but allow us going beyond a few cases to explore many interesting figures appeared in these rigid CY 3-folds. For the minimal number of intersecting points of a rational curve with the exceptional divisors in our framework, we obtain the complete solution on the location of $\PZ^1$ in $\widehat{T/G}$, and also the total number of those curves. In the case $G= \ZZ_3$, we have developed a quantitative  method of constructing  
rational curves in $\widehat{T/G}$, then carried it out in some cases to obtain an explicit form of the solution. 

The organization of this paper is as follows. 
In Section 2, we determine the complete structure of CY 3-folds $\widehat{T/G}$ for torus-orbifolds $T/G$ with only isolated singularities. There are only two such $\widehat{T/G}$: $T$ is either the triple-product of the 1-torus carrying an order 3 automorphism with $G= \ZZ_3$, or the Jacobian variety of Klein quartic with $G= \ZZ_7$, (see Examples 1 and 2 in paper). The derivation is based on two ingredients, one about the topological relation on the fixed-point set of an $n$-torus $T$ acted by $G$ with the $G$-invariant spaces of cohomologies of $T$, and the other on the classification of 3-tori carrying a cyclic automorphism group with a finite fixed-point set. The latter fact, which was contained in our previous work \cite{R7}, indeed can be understood as a consequence of CM type abelian varieties in \cite{ST} and the well-known fact on the ideal-class number of cyclotomic fields of small degree in number theory. We shall derive those facts relevant to the present paper in the appendix for easier reference. The crepant resolution $\widehat{T/G}$ of $T/G$ for both two cases are rigid CY 3-folds. For the next two sections, we study problems of rational curves in these two $\widehat{T/G}$ not contained in exceptional divisors. In Section 3, we formulate a precise problem of rational curves $C$ by imposing a condition on $C$ when intersecting the exceptional set of $\widehat{T/G}$, and introduce the counting number $d_C$ of exceptional components  meeting $C$. By considering the smooth model $\widehat{W_C}$ of the $G$-covering curve $W_C \subset T$ over the $C$-image in $T/G$, we investigate the general properties of $\widehat{W_C}$, in particular one obtains $d_C \geq 3$. For $d_C=3$, we are able to derive the smooth structure of $W_C$, the location of $C$ in $\widehat{T/G}$, and the total number of all such $C$'s.  In Section 4, we study the  existence problem of rational curves $C$ in $\widehat{T/G}$ with $d_C \geq 4$ in the case $G = \ZZ_3$, and the property $d_C \neq 4, 5$ is easily obtained. For $d_C=6$, $W_C$ is again a smooth curve, and we have constructed all $W_C$ in $T$, hence the rational curves $C$ in $\widehat{T/G}$. We also show that there is no $C$ with $d_C = 7$. Then we give the conclusion remarks in Section 5.

{\bf Notation.} To present our results, we prepare some notations. In this paper, $\ZZ, \RZ, \CZ$ will denote the ring of integers, real, complex numbers respectively, and we will write the order $m$ cyclic group $\ZZ/m\ZZ$ by $\ZZ_m$. For a topological space $X$, the Euler number of $X$ will be denoted by $\chi (X)$. 

\section{Rigid Calabi-Yau 3-folds}
It is known that the automorphism group of a polarized abelian variety is always finite. For the purpose of constructing CY space, we shall consider only those with automorphisms preserving the (holomorphic) volume form. 
  
Given a (compact) complex $n$-torus $T ( = \CZ^n/L )$ for $n \geq 2$, we denote by $Saut(T)$ the group of special automorphisms of $T$, i.e., the Lie-automorphisms preserving the volume $n$-form of $T$. A class of CY manifolds can be obtained by the following procedures: 
\bea(ll) 
(i) & {\rm Determine \ the \ structure \ of \ the \ orbifold}, \ T/G, \ {\rm for \ a \  (nontrivial) } \\
& {\rm  \ finite \ subgroup} \ G \subset Saut(T). \\
(ii) & {\rm Construct \ the \ (or \ a) \ crepant \ resolution} \ \widehat{T/G} \ {\rm of \ the \ orbifold} \ T/G.
\elea(T/G)      
For $n=2$, the well-known K3 surfaces are Kummer surfaces are obtained by the minimal resolution of $T$ quotiented by the involution. However for $n \geq 3$, the involution of a $n$-torus $T$ no longer serves the purpose of constructing CY $n$-folds as in K3 case, since it fails either the condition $(i)$ (for odd $n$) or the requirement $(ii)$ (for even $n$) in (\req(T/G)). Hence in searching CY $n$-fold by the method (\req(T/G)) for $n \geq 3$, the $n$-torus $T$ is required to possess certain extra special symmetries to generate the group $G$ in $(i)$, hence abelian varieties of CM-type naturally appear in this context. In general, the CM-structure of $T$ is related to the nature of  $T^g$ (the fixed point set of $g$) for $g \in G$, or equivalently, the $G$-fixed point set $T_G$,
$$
T_G : = \bigcup_{{\bf 1} \neq g \in G} T^g \ .
$$ 
A discussion of such relation was given in \cite{R7}. To make our presentation self-contained, we shall derive the relevant results in the appendix for easier use. 
For the construction of $(ii)$ in (\req(T/G)), the solution relies mainly on the solution of the local version, namely:  
\bea(ll)
(ii)'& {\rm Construct \ the \ (or \ a)  \ resolution} \ \widehat{\CZ^n/G} \ {\rm with \ trivial \ canonical \ bundle  } \\
& {\rm for \ a \ finite \ subgroup  } \ G \subset {\rm SL}_n ( \CZ ).
\elea(lii)
For $n \geq 4$, it is known that there is no crepant resolution $\widehat{\CZ^n/G}$ for a general subgroup $G$ of ${\rm SL}_n ( \CZ )$.  Hence the criterion of $G$ so that one could  obtain a positive result of (\req(lii)) becomes a non-trivial problem, of which little results have been known except a few special cases (e.g., see \cite{LR}). For $n=3$, the positive solution of (\req(lii)) has been known for all special group $G$ (see \cite{Rtop} and references therein). Hence, in principle one should be able to determine all the CY 3-folds arising from 3-torus-orbifolds through the procedure (\req(T/G)).      

First we describe certain topological constraints on a torus $T$ acted by a group $G$ such that $T/G$ has only isolated singularities. 
\begin{lemma} \label{lem:T_G}
Let $T$ be a $n$-torus, and $G$ be a finite (tori-)automorphism group of $T$, (no special condition on automorphisms required), such that $T_G$ consists of only finite elements. Then  

(i)The number $|T_G|$ is divisible by $|G|$ with $- \frac{|T_G|}{|G|}$ equal to the Euler number  
$\chi ((T \setminus T_G)/G )$.

(ii) The following equality holds:
$$
|T_G/ G| - \frac{|T_G|}{|G|} = 2 (1 + \sum_{j=2}^{n-1} {\rm dim} \ H^j (T)^G) + (-1)^n {\rm dim} \ H^n (T)^G  \ .
$$
(Here the suffix $G$ in the above right hand side means the $G$-invariant subspace.)
\end{lemma}
{\it Proof.} For an element $p \in T_G$, we consider a small $G$-invariant ball $B(p)$ in $T$ centered at $p$. As $G$ acts freely on $T \setminus T_G$, the equality, $|G| \chi ( (T \setminus T_G)/G ) = \chi ( T \setminus T_G )$, holds. Indeed, one has the following relation of the cohomology (with coefficients in $\CZ$):
\be
H^i ( (T \setminus T_G)/G ) = H^i ( T \setminus T_G  )^G \ \ \ \ {\rm for \ all} \ i \ .
\ele(TTG)
By the (cohomology) Mayer-Vietoris sequence of the pair $(T \setminus T_G , \cup_{p \in T_G} U(p) )$, one has 
$$
\chi ( T \setminus T_G ) + |T_G| = \chi (T) ~ (= 0 ) .
$$
The result $(i)$ then follows. The Mayer-Vietoris sequence of $((T \setminus T_G)/G , \cup_{p \in T_G} U(p)/G )$ gives the following relation:
$$
 \chi (T/G)  =  |T_G/ G| + \chi ( (T \setminus T_G)/G )  ~ (\stackrel{{\rm by} \ (i)}{=}  |T_G/ G| - \frac{|T_G|}{|G|}  ) \ ,
$$
By comparing the Mayer-Vietoris sequences of $(T \setminus T_G , \cup_{p \in T_G} U(p) )$ and $((T \setminus T_G)/G , \cup_{p \in T_G} U(p)/G )$, then using (\req(TTG)), we have the identification of vector spaces:
$$
H^i (T/G)  = H^i(T)^G  \ ~ ~ {\rm for} \ 0 \leq i \leq 2n \ .
$$
Using $H^j(T)^G \simeq H^{2n-j}(T)^G$ for $0 \leq j \leq n$, and $H^1(T)^G = 0$ ( by the finiteness of $T_G$), one can express $\chi (T/G)$ by
$$
\chi (T/G) = 2 (1 + \sum_{j=2}^{n-1} {\rm dim} \ H^j (T)^G) + (-1)^n {\rm dim} \ H^n (T)^G \ .
$$
Then $(ii)$ follows. 
$\Box$ \par \vspace{.2in} 
For the rest of this section, $T$ will always denote a 3-torus, and $G$ a finite special subgroup of $T$ with $|T_G| < \infty$. One may identify the universal cover $\widetilde{T}$ of $T$ with the tangent space of $T$ at the identity element $o$ of $T$. Then  the group $G$ can be regarded as a subgroup of ${\rm SL}_3(\CZ)$ by the differential at $o$:
\be
G \hookrightarrow {\rm SL}(\widetilde{T}) ~ (\simeq \ {\rm SL}_3(\CZ)) \ , \ \ g \mapsto (dg)_o \ .
\ele(do)
For $p \in T_G$, the isotropy subgroup of $G$ at $p$ will be denote by $G_p$.

Claim: The order of every non-trivial element $g \in G$ must be odd, and  ${\rm dim} \ H^3 (T)^G  = 2$. In fact, if the order of $g$ is an even number $2k (>1)$, then $g^k$ is an order 2 element in $G$. This implies the differential $d (g^k)_o$ of $g^k$ in (\req(do)) with the eigenvalues, $-1, -1, 1$, hence $g^k$ fixes a 1-subtorus of $T$, a contradiction to $|T_G| < \infty$. By $G \subset Saut(T)$, $H^3 (T)^G$ contains the 2-dimensional subspace $H^{3,0}(T) +  H^{0,3}(T)$. In $G$, we consider a non-trivial element $g$, with the representation in (\req(do)) given by a diagonal transformation for some coordinates $(z_1, z_2, z_3)$ of $\CZ^3$. Among all the 3-forms of $T$, $dz_1 \wedge dz_2 \wedge dz_3$ and $\overline{dz}_1 \wedge \overline{dz}_2 \wedge \overline{dz}_3$ are the only $g$-invariant eigen-forms. Otherwise, the diagonal action of $g$ on $\CZ^3$ has $-1$ as one of its eigenvalues, which implies $g^2$ is a non-trivial automorphism fixing a 1-subtorus of $T$, again a contradiction to $|T_G| < \infty$. Therefore,  $H^3(T)^g$ is a 2-dimensional subspace, so is $H^3 (T)^G $.  By which, the equality in Lemma \ref{lem:T_G} $(ii)$ becomes  
\be
|T_G/ G| - \frac{|T_G|}{|G|} = 2 \ {\rm dim} \ H^2 (T)^G    \ .
\ele(TGG3)

We now describe two standard examples of 3-tori (which appeared previously in \cite{RY}) for later use. 

{\bf Example 1}. $E (\omega)^3$ where $\omega := e^{\frac{2 \pi {\rm i}}{3}}$. Let $E(\omega)= \CZ/ (\ZZ + \ZZ \omega)$ be the 1-torus with the order 3 automorphism $m_{\omega}$ (defined by $[z] \mapsto [\omega z]$), and $E (\omega)^3$ the triple-product of $E (\omega)$. The triple-product of $m_\omega$ defines the automorphism of $E (\omega)^3$,
$$
m_{\omega}^3 : E (\omega)^3 \longrightarrow E (\omega)^3 \ , \ \ ([z_1], [z_2], [z_3] ) \mapsto ([\omega z_1], [\omega z_2], [\omega z_3] ) \ ,
$$  
which generates an order 3 special automorphism group $G = \langle m_\omega^3 \rangle$ of $E (\omega)^3$. Then the map (\req(do)) is given by $(d m_{\omega}^3)_0 = \omega I_3$, hence $H^2 (T)^G $ is the 9-dimensional space with a basis given by
$$ 
H^2 (T)^G = \bigoplus_{i, j=1}^3 \CZ dz_i \wedge \overline{dz}_j \ .
$$
The fixed point set $(E (\omega)^3)_G$  is the triple-product of $E(\omega)_{ \langle m_{\omega} \rangle }$ (= the  order 3 additive subgroup  of $E (\omega)$ generated by $[\frac{1+2 \omega}{3}]$). For each $p \in (E (\omega)^3)_G$, we have $G_p = G$, which acts as the $\omega$-multiplication near $p$ in $T$. By which, one can also verify the formula (\req(TGG3)) directly from  the terms involved. 
$\Box$ \par \vspace{.2in}
{\bf Example 2}. $A( \QZ (\mu ))$ where $\mu := e^{\frac{2 \pi {\rm i}}{7}}$. Let $\varphi_i$, $1 \leq i \leq 3$, be the three field-embeddings of $\QZ(\mu)$ into $\CZ$ with $\varphi_1 (\mu) = \mu$, $\varphi_2 (\mu) = \mu^2$ and $\varphi_3 (\mu) = \mu^4$. Then $(\QZ(\mu) , \{ \varphi_i\}_{i=1}^3)$ is a CM-field \cite{Sh}. The $\RZ$-isomorphism, $u : \QZ(\mu) \otimes \RZ \longrightarrow \CZ^3$, which embeds $\QZ(\mu)$ into $\CZ^3$ via  $u (x) = (\varphi_1(x) , \varphi_2 (x), \varphi_3 (x) )^t$ for $x \in \QZ(\mu)$, gives rise to a lattice $ u( \ZZ[\mu] )$ in $\CZ^3$ with the integral basis $\{u( \mu^i) \}_{i=0}^5$ expressed by
$$
\bigg( u(1), u(\mu), \cdots , u(\mu^5) \bigg)= \left( \begin{array}{cccccc}
1 & \mu & \mu^2 & \mu^3 & \mu^4 & \mu^5 \\
1 & \mu^2 & \mu^4 & \mu^6 & \mu & \mu^3 \\
1 & \mu^4 & \mu & \mu^5 & \mu^2 & \mu^6 \\
\end{array} \right) \ .
$$
The 3-torus
$$
A( \QZ (\mu )) = \CZ^3 / u( \ZZ[\mu] )
$$
is an abelian variety of CM-type with a special automorphism group $G$ generated by
$$
m_{\mu} : A( \QZ (\mu )) \longrightarrow A( \QZ (\mu )) \ , \ \ [(z_1, z_2, z_3)^t] \mapsto
[(\mu z_1, \mu^2 z_2, \mu^4 z_3)^t] \ .
$$ 
Then the fixed point set $A( \QZ (\mu ))_G$ is an order 7 subgroup of $A( \QZ (\mu ))$ with the generator $\frac{1}{7}\sum_{j=0}^5 (j+1) u(\mu^j)$. For each $p \in A( \QZ (\mu ))_G$, one has $G_p = G$, acting on $T$ near $p$ as the diagonal matrix with eigenvalues $\mu, \mu^2 , \mu^4$. Then $H^2 (T)^G$ is a 3-dimensional space with a basis $\{ dz_j \wedge \overline{dz}_j \}_{j=1}^3$. By which, one can also directly check the relation (\req(TGG3)).
$\Box$ \par \vspace{.2in}
We now show that the above two examples are indeed the only ones which give rise to CY 3-folds by resolving singularities of the torus-orbifold.  
\begin{theorem} \label{thm:RCY}
Let $T$ be a 3-torus and $G$ a finite subgroup of $ Saut(V)$ with $|T_G| < \infty$. Then

(i) $G$ is a cyclic group of order $3$ or $7$. The pair $(T, G)$ is isomorphic to either $(E (\omega)^3, \langle m_\omega^3 \rangle)$ (in {\rm Example 1}), or $(A( \QZ (\mu )), \langle m_{\mu}  \rangle)$ (in {\rm Example 2}).

(ii) There exists the unique crepant (toric) resolution $\widehat{T/G}$ of $T/G$. Both $\widehat{T/G}$ are rigid projective CY 3-folds with the zero Hodge numbers except $h^{0,0} = h^{3,3}= h^{3,0}= h^{0,3} =1$, and the following ones:
$$
\begin{array}{lll}
\widehat{E (\omega)^3/ \langle m_\omega^3 \rangle} : & h^{1,1} = h^{2,2}= 36 & ( \chi = 72) \ ; \\
\widehat{A( \QZ (\mu ))/ \langle m_{\mu}  \rangle} : & h^{1,1} = h^{2,2}= 24 & ( \chi = 48) \ .
\end{array}
$$
\end{theorem}
{\it Proof}. $(i)$. First we consider the case when $G$ is a cyclic group. By using results on abelian varieties with complex multiplication \cite{ST}, one can classify all the 3-tori $T$ with a finite automorphism group $G$ of order $\geq 3$  \cite{R7}. By which, $(E (\omega)^3, \langle m_\omega^3 \rangle)$ and $(A( \QZ (\mu )), \langle m_{\mu}  \rangle)$ are the only $(T, G)$ with $G \subset Saut(T)$ and $|T_G| < \infty$. (We shall derive the relevant facts in the appendix for the sake of completeness of this paper). Hence it suffices to show that there is no non-cyclic special automorphism group $G$ for a 3-torus $T$ with $|T_G| < \infty$. Otherwise, let $({\cal T, G})$ be a such pair with the minimal order $|{\cal G}|$ among all $(T, G)$ such that $G$ is a non-cyclic special automorphism group of $T$ with $|T_G| < \infty$.
By the finiteness of ${\cal T_G}$, an order 3 subgroup of ${\cal G}$, if exists, must be generated by the multiplication of $\omega$ via the embedding (\req(do)), hence unique and contained in the center of ${\cal G}$. This implies that every element in ${\cal G}$ must be of order 7 from our previously known classification of 3-tori with a cyclic special automorphism group.  
Therefore $|{\cal G}| = 7^k $ for some $k \geq 2$, and ${\rm dim} \ H^2 (T)^G  = 3$. By the minimal property of $|{\cal G}|$, the isotropy subgroup ${\cal G}_p$ for an element $p \in {\cal T_G}$ is either equal to ${\cal G}$, or a subgroup of order 7. Denote by ${\cal T_G}_0$, ${\cal T_G}_1$ the subsets of ${\cal T_G}$ consisting of the elements $p$ with $|{\cal G}_p|= 7^k$ or 7 respectively. ( Note that $ o \in {\cal T_G}_0$.)  Then $|{\cal T_G}_1|$ is divisible by $7^{k-1}$, and $|{\cal T_G}/{\cal G}| = |{\cal T_G}_0| + \frac{|{\cal T_G}_1|}{7^{k-1}}$. By Theorem \ref{lem:T_G} $(i)$, $|{\cal T_G}| (= |{\cal T_G}_0| + |{\cal T_G}_1|)$ is divisible by $7^k$, which implies $\frac{|{\cal T_G}_0|}{7^{k-1}}$ is a positive integer.
Using (\req(TGG3)) for the pair $({\cal T , G})$, one has 
$$
(\frac{7^{k-1}-1}{2}) \frac{|{\cal T_G}_0|}{7^{k-1}} + 3 \frac{|{\cal T_G}_0| + |{\cal T_G}_1|}{7^k } = 3    \ ,
$$
which leads to a contradiction as all terms in the above left hand side are positive integers. The result $(i)$ then follows.

$(ii)$ By the description of $T_G$ in Examples 1 and 2, $T/G$ is an orbifold with only isolated singularities, and the local structures near singular points are all the same, given by $\CZ^3/G$ via (\req(do)). We have $|{\rm Sing}(T/G)| =27, G= \langle \omega I \rangle$ in Example 1, and $|{\rm Sing}(T/G)| =7, G= \langle {\rm dia} [ \mu, \mu^2, \mu^4 ] \rangle$ in Example 2. For both cases, one has the crepant toric resolution $\widehat{\CZ^3/G}$ of $\CZ^3/G$ ( see, e.g., \cite{R89} Proposition 2.), by which one obtains the CY resolution $\widehat{T/G}$ over $T/G$,
$$
\pi: \widehat{T/G} \longrightarrow T/G \ .
$$
Indeed, the crepant toric resolution $\widehat{\CZ^3/G}$ is unique in both cases. The exceptional set $E$ in $\widehat{\CZ^3/G}$ is given by $E = \PZ^2$ in Example 1, and $E = D_1 + D_2 + D_3$ in Example 2, where $D_i \simeq \overline{\HZ^{-2}}$, the $\PZ^1$-bundle over $\PZ^1$ by compactifying the $(-2)$-hyperplane bundle with $\infty$-section, are normal crossing divisors with the only intersections as follows:
$$
\begin{array}{lll}
D_1 \cdot D_2 &=  0{\rm -section \ of} \ D_1 &= {\rm a \ \PZ^1-fiber \ of } \ D_2 , \\
D_2 \cdot D_3 &=  0{\rm -section \ of} \ D_2 &= {\rm a \ \PZ^1-fiber \ of } \ D_3 , \\
D_3 \cdot D_1 &=  0{\rm -section \ of} \ D_3 &= {\rm a \ \PZ^1-fiber \ of } \ D_1 ,
\end{array}
$$
(for the derivation, see, e.g.  \cite{RY}  \S 3 Example II).  Hence in both case, we have
$$
H^1 ( \widehat{\CZ^3/G} ) = H^3 ( \widehat{\CZ^3/G} ) = 0 \ , \ \ \ H^2 (\widehat{\CZ^3/G}) = \sum_{D} \CZ c(D) \ 
$$
where the index $D$ runs through all irreducible exceptional divisors, and $c(D)$ denotes the first-Chern class of the line bundle ${\cal O}[D]$. 
Choose a sufficient small $G$-invariant ball $U(p)$ in $T$ near each $p \in T_G$ as before (note that $G_p =G$ in the present situation),  and denote the $G$-quotient of $U(p)$ in $T/G$ around $[p] (=G$-orbit of $p$) by $V([p])= U(p)/G$. Note that one can make the identification: $\widehat{T/G} \setminus \pi^{-1}({\rm Sing}(T/G))= (T \setminus T_G)/G$, and $H^j ((T \setminus T_G)/G) = H^j (T)^G$ for $j=1,2,3$. By Mayer-Vietoris sequence of the pair $\bigg(\widehat{T/G} \setminus \pi^{-1}({\rm Sing}(T/G)) , ~ \cup_{[p] \in {\rm Sing}(T/G)} \pi^{-1}(V([p]) \bigg)$, one obtains
$$
\begin{array}{ll}
H^1 (\widehat{T/G})  & =  H^1 (T)^G \ ~ ( = 0 ) , \\ 
H^2 (\widehat{T/G})  & = H^2(T)^G + \sum_{D} \CZ c(D) \ , \ \ {\rm where} \ D: {\rm irred. \ exceptional \ divisor \ in } \ \widehat{T/G} \ , \\
H^3 (\widehat{T/G})  & =  H^3 (T)^G \ ~ ( = \CZ v_T + \CZ \overline{v_T} \ , \ \ v_T : {\rm holomorphic \ volume \ form \ of \ } T ) .
\end{array}
$$
Since elements in ${\rm Sing} (T/G)$ and $T_G$ are in one-one correspondence, by the descriptions of $H^2(T)^G$ and $T_G$ in Examples 1 and 2, plus the local structure of exceptional divisors in the crepant resolution $\widehat{T/G}$ over ${\rm Sing} (T/G)$, the Hodge numbers of $\widehat{T/G}$ in  $(ii)$ then follow immediately.
$\Box$ \par \vspace{.2in} \noindent
{\bf Remark}. By the structure of $(T, G)$ in Theorem \ref{thm:RCY} $(i)$, translations by elements in $T_G$ and the centralizer of $G$ in $Saut(T)$ give rise to biregular automorphism groups of $T/G$, hence ones on the CY 3-fold $\widehat{T/G}$. For the case in Example 1,   
$\widehat{E (\omega)^3/ \langle m_\omega^3 \rangle}$ possesses the infinite automorphism group
$(E(\omega)^3)_G \bullet {\rm PSL}_3(\ZZ[\omega])$ (=
the semidirect product of ${\rm PSL}_3(\ZZ[\omega]$ by $(E(\omega)^3)_G$)
$\Box$ \par \vspace{.2in} 
There are several representations for the isomorphic class of $( A( \QZ (\mu )) , \langle m_\mu \rangle )$ in Example 2. First, we show an easy lemma for later use:
\begin{lemma} \label{lem:g3G7}
Let $W$ be a genus $3$ curve with an order $7$ automorphism group ${\bf g}$. Then $W$ is birational to the affine curve $y^7 = x^k (x-1)$ in $\CZ^2$ for $k=1, 2, 3, 5$, and ${\bf g}$
is generated by the morphism $(x, y) \mapsto (x , \mu y)$. Furthermore, a basis of abelian differentials of first kind of $W$ is given by the following list:
$$
\begin{array}{| c |c | }
\hline 
k & {\rm forms \ for} \ y^7 = x^k (x-1) \\
\hline 
1 & dx / y^4 , ~ ~ dx/ y^5 , ~ ~ dx/ y^6 \\
\hline 
2 & dx / y^3 , ~ ~ xdx/ y^5 , ~ ~ xdx /y^6  \\
\hline 
3 & dx / y^2 , ~ ~ xdx / y^4 , ~ ~ x^2dx / y^6 \\
\hline 
5 & x^2dx / y^4 , ~ ~ x^3dx / y^5  , ~ ~ x^4dx /y^6 \\
\hline 
\end{array}
$$
As a consequence, $y^7 = x^2 (x-1)$ is the only case such that ${\bf g}$
induces the special automorphism group of the Jacobian of $W$.
\end{lemma}
{\it Proof.} By Hurwitz's theorem, one has the Galois ${\bf g}$-cover of $W$ over $W/{\bf g} = \PZ^1$ with three branched points, which can be assumed to be $0, 1, \infty$. Then $W$ is birational to $y^7 = x^k (x-1)$ for a positive integer $k < 7$, and the morphism, $(x, y) \mapsto (x , \mu y)$, gives rise to a generator of ${\bf g}$. As $\infty$ is a branched point of the cover $W$ over $\PZ^1$, $k$ can not be $6$. Furthermore, $y^7 = x^k (x-1)$ for $k=2$ and 4 are birational curves via the birational map: $(x , y ) \mapsto (\frac{1}{x}, \frac{-\mu y}{x})$. The basis of abelian differentials of first kind in the above table can be obtained by the direct computation on each case.
$\Box$ \par \vspace{.2in} \noindent
By the above lemma, the class $( A( \QZ (\mu )) , \langle m_\mu \rangle )$ can be represented by the Jacobian of the curve $y^7 = x^2 (x-1)$. One can also describe this curve by the Klein quartic \cite{Kl}:
$$
K : \ \ Z_1^3 Z_2 + Z_2^3 Z_3 + Z_3^3 Z_1 = 0 \ , \ \ [ Z_1 , Z_2 , Z_3 ] \in \PZ^2 \ .
$$
It is known that $K$ is the genus 3 curve having the automorphism group $Aut(K)$ isomorphic to the order 168 simple group. Hence the Jacobian ${\rm Jac}(K)$ of Klein quartic carries the special automorphism group with the same structure as $Aut(K)$. Inside $Aut(K)$, there is an order 7 subgroup, denoted by $\kappa_7$. By Theorem \ref{thm:RCY}, we have 
$$
( A( \QZ (\mu )) , \langle m_\mu \rangle ) ~ \simeq ({\rm Jac}(K) , \ \kappa_7 ) \ , \ \ \ \widehat{A( \QZ (\mu ))/\langle m_\mu \rangle} ~ \simeq \widehat{{\rm Jac} (K) / \kappa_7}   \ .
$$
By Lemma \ref{lem:g3G7}, $K$ can be represented by the affine curve $y^7 = x^2 ( x - 1 )$, with $\kappa_7$ corresponding to the cyclic group generated by the map, $(x, y) \mapsto (x, \mu y)$.  Another representation of the torus in Example 2 is the 3-product of elliptic curve $E( \eta ) = \CZ/(\ZZ+\ZZ \eta)$ where $\eta := \frac{-1 + \sqrt{7}{\rm i}}{2} (= \mu + \mu^2 + \mu^4)$. As an element $\beta$ in $\ZZ[\eta]$ gives rise to an endomorphism of $E( \eta )$ induced by the $\beta$-multiplication on $\CZ$, one has $End (E(\eta) ) = \ZZ[\eta]$. Hence $Saut ( E(\eta)^3 ) = {\rm SL}_3 (\ZZ[\eta])$. Consider the order 7 special automorphism group $G_{\eta}$ of $E(\eta)^3$ generated by the following matrix in ${\rm SL}_3 (\ZZ[\eta])$,
$$
 \left( \begin{array}{ccc} 0 & 0 & 1 \\
1 & 0 & \eta + 1 \\
0 & 1 & \eta \end{array}
\right) \ .
$$
Then by Theorem \ref{thm:RCY}, we have 
$$
( A( \QZ (\mu )) , \langle m_\mu \rangle ) ~ \simeq (E(\eta)^3  , \ G_{\eta} ) \ , \ \ \ \widehat{A( \QZ (\mu ))/\langle m_\mu \rangle} ~ \simeq \widehat{E(\eta)^3/G_{\eta}}  \ .
$$

\section{Rational Curves in $\widehat{T/G}$}
For the rest of this paper, we are going to study a rational curve problem in the rigid CY 3-folds $\widehat{T/G}$ in Theorem \ref{thm:RCY}. As all the exceptional divisors in $\widehat{T/G}$ are smooth rational surfaces, we shall discuss those rational curves in $\widehat{T/G}$ not contained in the exceptional set. Indeed, we consider only the rational curves $C$ in $\widehat{T/G}$ with the following property:\par \vspace{.1in} \noindent
$(P)$: $C$ is a $\PZ^1$-curve intersecting each connected component of the exceptional set in $\widehat{T/G}$ at most at one point.  
\par \vspace{.08in} \noindent
For the rest of this paper, we shall denote by $C$ a $\PZ^1$-curve in $\widehat{T/G}$ with the property $(P)$.

Denote the map of $G$-orbits in $T$ by
$$
\wp : T \longrightarrow T/G \ .
$$
For a $\PZ^1$-curve $C$  in $\widehat{T/G}$ with $(P)$, its image $\pi (C)$ in $T/G$ is a rational curve, but possibly singular when passing through ${\rm Sing} (T/G)$. We shall denote 
$$
\begin{array}{ll}
d_C :&= |\pi (C) \cap {\rm Sing} (T/G)| \\
&= | \{ {\rm connected \ components \ in \ the \ exceptional \ set \ intersecting } \ C \}|  \ ~ , ~ ( {\rm by} \ (P) )   \ ,
\end{array}
$$
and $W_C := \wp^{-1}( \pi (C))$, which is a curve in the torus $T$.  
\begin{lemma} \label{lem:3pts}
$W_C$ is irreducible, and  $d_C \geq 3$.
\end{lemma}
{\it Proof.} If $W_C$ is reducible,  it must consist of $|G|$ irreducible components as the group $G$ acting on $W_C$ is of order 3 or 7. Then each irreducible component is a (possibly singular) rational curve, and can be lifted to the universal cover $ ( \simeq \CZ^3)$, which is impossible by the fact that there is no non-constant global function for compact Riemann surfaces. Therefore $W_C$ is irreducible. Note that elements of $T_G$ and ${\rm Sing}(T/G)$ are in one-one correspondence under $\wp$. If $d_C \leq 1$, $W_C$ must consists of $|G|$ irreducible components by the simply-connectness of $\PZ^1$ and $\CZ$; while in the case $d_C=2$, $W_C$ is again a rational curve (possibly with at most two singularities), hence can be lifted to  the universal cover of $T$. Both situations are impossible. Hence one obtains $d_C \geq 3$.
$\Box$ \par \vspace{.2in} \noindent
By the relation between $T_G$ and ${\rm Sing}(T/G)$, we have $|W_C \cap T_G| = d_C$, and ${\rm Sing} (W_C) \subset W_C \cap T_G $. The normalization of $W_C$, denoted by $\widehat{W_C}$, is a smooth curve with a bijective morphism onto $W_C$. The similar relation also holds between $C$ and $\pi (C)$, and one has the Galois  $G$-cover of $\widehat{W_C}$ over $C$ with $d_C$ branched elements corresponding those in $W_C \cap T_G$,
\be
\tau_C :  \widehat{W_C} \longrightarrow C \ .
\ele(tauC)
Hence $\widehat{W_C}$ is a Riemann surface of genus $d_C -2$ for $|G|=3$, and $3d_C -6$ for $|G|=7$.   

We now discuss the structure of the genus $(d_C -2)$ Riemann surface $\widehat{W_C}$ in the case of Example 1 for $d_C \geq 3$. Denote by $\theta$ the order 3 automorphism of $\widehat{W_C}$ induced by $m_\omega^3$ of $E(\omega)^3$. One has the following local description of the morphism from $\widehat{W_C}$ into $E(\omega)^3$ near $(\widehat{W_C})_{\langle \theta \rangle}$.
\begin{lemma} \label{lem:lf3}
Let $s$ be an element in $(\widehat{W_C})_{\langle \theta \rangle}$, and $(s_1, s_2, s_3) \in E(\omega)^3$ be its corresponding element in $W_C$. Let $z_j $ be the uniformizing coordinate of $E(\omega)$ centered at $s_j$ for $j=1,2,3$. Then there exists a local coordinate $t$ in $\widehat{W_C}$ centered at $s$ with the following local expression of $\widehat{W_C}$ into $E(\omega)^3$ near $s$:
$$
t \mapsto (z_1, z_2 , z_3) = ( \alpha_1 t^{r_{s,1}+1} f_1(t^3) , ~ \alpha_2 t^{r_{s,2}+1} f_2(t^3), ~\alpha_3 t^{r_{s,3}+1} f_3(t^3) )  \ , 
$$ 
for some $\alpha_j \in \CZ $, $ r_{s,j} \in \ZZ_{\geq 0}$, and analytic functions $f_j(*)$  with $f_j(0)=1$, such that either all $r_{s,j} \equiv 0 \pmod{3}$ or all $r_{s,j} \equiv 1 \pmod{3}$. (Note that the condition of $r_{s, j}$'s depends only on those $j$ with $\alpha_j \neq 0$.) Furthermore, one has 
\be
{\rm gcd}. \{ r_{s,j} + 1 \ | 1 \leq j \leq 3 \ , \ \ \alpha_j \neq 0 \} = 1 \ . 
\ele(gcdr)
In particular, if $(s_1, s_2, s_3)$ is non-singular in $W_C$, one has $r_{s,j} \equiv 0 \pmod{3}$ for all $j$. 
\end{lemma}
{\it Proof}. Choose a local coordinate $t$ in $\widehat{W_C}$ centered at $s$ so that $\theta$ is locally described by $t \mapsto \omega t$, or $t \mapsto \omega^2 t$. Then the local expression of $\widehat{W_C}$ into $E(\omega)^3$ near $s$ is given by  
$$
t \mapsto (z_1, z_2 , z_3) = ( \alpha_1 t^{r_{s,1}+1}F_1(t) , ~ \alpha_2 t^{r_{s,2}+1} F_2(t), ~ \alpha_3 t^{r_{s,3}+1}F_3(t) )
$$
for some $\alpha_j \in \CZ $, $ r_{s,j} \in \ZZ_{\geq 0}$, and analytic functions $F_j(t)$  with $F_j(0)=1$. Since the above map sends $\theta (t)$ to $(\omega z_1, \omega z_2, \omega z_3 )$,  $F_j(t)$ must depend only on $t^3$, and  the integers $r_{s,j}$ are all $ \equiv 0 \pmod{3}$ or all $\equiv 1 \pmod{3}$ according to $\theta(t)= \omega t$ or $ \omega^2 t$ respectively. The condition (\req(gcdr)) follows from the bijective relation between $\widehat{W_C}$ and $W_C$. When $W_C$ is non-singular at the element $(s_1, s_2, s_3)$, there exists some $j$ with $\alpha_j \neq 0$ and $r_{s,j}=0$. The result then follows. 
$\Box$ \par \vspace{.2in} 
The projection of $W_C$ to each  factor of $E (\omega)^3$ gives rise to the $\ZZ_3$-equivariant morphism,
\be
\rho_j : \widehat{W_C} \longrightarrow E ( \omega ) \ , \ \ ~ \ j = 1, 2, 3 ,
\ele(rhoj)
by identifying the group $\ZZ_3$ with  $\langle \theta \rangle$ on $\widehat{W_C}$, and $\langle m_\omega \rangle$ on $E(\omega)$. The curve $W_C$ is described by all $\rho_j$'s, and $W_C/\ZZ_3$  gives rise to the $C$ in $\widehat{T/G}$ for $T= E(\omega)^3$.
Consider those non-constant $\rho_j$s, and denote a such morphism $\rho_j$ simply by $\rho$ for the sake of convenient notations.
$$
\rho: \widehat{W_C} \longrightarrow E ( \omega ) \ .
$$
The $\ZZ_3$-quotient of $\rho$ gives rise to a rational function of $C \ (\simeq \PZ^1) $: 
$$
\overline{\rho}: C \longrightarrow \PZ^1 \ .  
$$
Note that $|(\widehat{W_C})_{\ZZ_3}|= d_C $ and $|E ( \omega )_{\ZZ_3}|=3$. 
As $(\widehat{W_C})_{\ZZ_3}$ is sent to $E ( \omega )_{\ZZ_3}$ by $\rho$, 
the group $\ZZ_3$ acts freely on $\rho^{-1}(E ( \omega ) \setminus E ( \omega )_{\ZZ_3})$ and $E ( \omega ) \setminus E ( \omega )_{\ZZ_3}$. Hence the morphisms $\rho$ and $\overline{\rho}$ have the same degree, denoted by 
$$
\delta = {\rm deg} (\rho) = {\rm deg} (\overline{\rho}) \ . 
$$ 
The  $\rho$-ramification index at $s \in (\widehat{W_C})_{\ZZ_3}$ is denoted by  $r_s +1$ for some $r_s \geq 0$. (Note that $r_s$ is the $r_{s,j}$ in Lemma \ref{lem:lf3}.) 
\begin{lemma} \label{lem:rhoI}
We have $\rho^{-1} ( E ( \omega )_{\ZZ_3} ) = (\widehat{W_C})_{\ZZ_3}$. For each $e \in E ( \omega )_{\ZZ_3}$, the degree $\delta$ of $\rho$ (or $\overline{\rho}$) can be expressed by the ramification indices of elements in $\rho^{-1}(e)$,  
$$
\delta = \sum_{s \in \rho^{-1}(e)} (r_s+1) \ .
$$  
\end{lemma}
{\it Proof}. For $e \in E ( \omega )_{\ZZ_3}$, we denote $O_e := \rho^{-1}(e) \setminus (\widehat{W_C})_{\ZZ_3}$, and the $\rho$-ramification index at $x \in O_e$ by $r_x+1$ for some $r_x \geq 0$. By the structure of  $\rho$-fiber over a generic point near $s$,  the degree of $\rho$ is expressed by 
$$
\delta = \sum_{s \in \rho^{-1}(e) \cap (\widehat{W_C})_{\ZZ_3}} (r_s+1) + \sum_{x \in O_e} (r_x + 1) \ .
$$
Under the $\ZZ_3$-quotient, elements of $E ( \omega )_{\ZZ_3}$ and $E ( \omega )_{\ZZ_3}/\ZZ_3$ are in one-one correspondence, and the same for    $(\widehat{W_C})_{\ZZ_3}$ and $(\widehat{W_C})_{\ZZ_3}/ \ZZ_3 ( \subset C )$. Furthermore, $\rho$ and $\overline{\rho}$ have the same ramification indices at the corresponding $\ZZ_3$-fixed elements.
The group $\ZZ_3$ acts freely on $O_e$ which consists of $\ZZ_3$-orbits, and all elements in a $\ZZ_3$-orbit have the same ramification index. Hence $\sum_{x \in O_e} (r_x + 1)$ is divisible by 3, and the degree of $\overline{\rho}$ is given by
$$
\delta = \sum_{s \in \rho^{-1}(e) \cap (\widehat{W_C})_{\ZZ_3}} (r_s+1) + \frac{\sum_{x \in O_e} (r_x + 1)}{3} \ .
$$
By comparing the above two expressions of $\delta$, one concludes that $O_e$ is an empty set. Then the results easily follow.
$\Box$ \par \vspace{.2in} \noindent
Denote by ${\cal B}$ the collection of all branched elements in $\widehat{W_C} \setminus (\widehat{W_C})_{\ZZ_3}$. For $q \in {\cal B}$, the $\rho$-ramification index at $q$ is expressed by $r_q + 1$ for some $r_q \geq 1$. As $\ZZ_3$ acts freely on ${\cal B}$, there exist some non-negative integers $\gamma^{(0)}$ and $\gamma^{(1)}$ so that 
$$
\gamma^{(0)} = \sum_{s \in (\widehat{W_C})_{\ZZ_3}} r_s   ,  \ \ \ ~ ~
3 \gamma^{(1)} = \sum_{p \in {\cal B}} r_p       \ \ .  
$$ 
For the morphism $\overline{\rho}$, all the branched elements in $C$ are descended from those in $\widehat{W_C}$ through $\tau_C$ ( in (\req(tauC))).  
We have 
$$
{\rm the } \ \overline{\rho}{\rm -ramification \ index \ at} \ \tau_C(p) = {\rm the } \  \rho {\rm -ramification \ index  \ at } \  p \ , \ \ {\rm for } \ p \in (\widehat{W_C})_{\ZZ_3} \cup {\cal B} \ . 
$$ 
Then by using the Hurwitz's formula for $\rho$ and $\overline{\rho}$, one obtains
\be 
 2 d_C - 6 =   \gamma^{(0)} + 3 \gamma^{(1)}   , \ \ \ \ \
2 \delta -2 =  \gamma^{(0)} + \gamma^{(1)}   \ , \ ~ ~ ({\rm hence} \ \delta = d_C -2 - \gamma^{(1)}) \ .
\ele(grho)
\begin{lemma} \label{lem:Ed345}
Let $\rho_j$s be the morphisms in $(\req(rhoj))$.

(i) When $d_C = 3$, $\rho_j$ is either  a constant-map with the value in $E(\omega)_{\langle m_\omega \rangle }$, or a biregular morphism. 

(ii) $d_C \neq 4 , 5$. 

(iii) If $W_C$ is a smooth curve in $T$,  $d_C$ is divisible by 3. 
\end{lemma}
{\it Proof.} Note that a constant $\rho_j$ must take the value in $E(\omega)_{\langle m_\omega \rangle }$. When $d_C = 3$, one has $\delta=1$ by (\req(grho)) for a non-constant $\rho_j (= \rho)$. The result $(i)$ then follows. 

$(ii)$. When $d_C \geq 4$, $\widehat{W_C}$ has the genus greater than 1. There exist at least two non-constant $\rho_j$s. If $d_C=4$, by the relation (\req(grho))  we have  $\gamma^{(0)} = \delta = 2$ and $|{\cal B}| = 0$ for every non-constant morphism $\rho_j$, hence with only two ramified elements, all belonging to $(\widehat{W_C})_{\ZZ_3}$. Say $\rho_1$ to be a non-constant morphism. Let $s$ be an element in $(\widehat{W_C})_{\ZZ_3}$ with the $\rho_1$-ramification  index 2. By the property (\req(gcdr)), either $\rho_2$ or $\rho_3$ has the ramification index at $s$ equal to one. This leads to a contradiction by the modular-3-property for all $\rho_j$-ramification indices at $s$ in Lemma \ref{lem:lf3}. 
 
If $d_C=5$, the relation (\req(grho)) has only two solutions : 
$$
\begin{array}{ll}
\gamma^{(0)} = \gamma^{(1)}  = 1  \ , \ \delta =2 \ ; & {\rm or} \ \ \
\gamma^{(0)}=4 \ , \ |{\cal B}|=  0 , \  \delta =3 \ .
\end{array}
$$
In the case $\gamma^{(0)}=4$, by $\delta=3$ there is an element in $(\widehat{W_C})_{\ZZ_3}$ with the ramification index 3, which is impossible by Lemma \ref{lem:lf3}. Therefore all non-constant $\rho_j$ are described by the first solution with $\delta =2$ . By (\req(gcdr)), for each $s \in (\widehat{W_C})_{\ZZ_3}$ there exists one $\rho_j$ with the ramification index 1 at $s$. This contradicts the condition $\gamma^{(0)}=1$ by the "the modular-3-property of ramification indices" in Lemma \ref{lem:lf3}. The result $(ii)$ then  follows.

$(iii)$. When $W_C$ is non-singular, all $\rho_j$-ramification indices at each $s \in (\widehat{W_C})_{\ZZ_3}$ are of the form $r_s+1$ with $r_s \equiv 0 \pmod{3}$ by Lemma \ref{lem:lf3}. Hence $\gamma^{(0)}$ is divisible by 3 in the first relation of (\req(grho)), which implies $d_C$ is divisible by 3.    
$\Box$ \par \vspace{.2in} 

For $(T, G)$ in Example 2, we have the following lemma.
\begin{lemma} \label{lem:Tmu}
For $(T, G)= A( \QZ (\mu )) , \langle m_\mu \rangle ) $, let $C$ be a rational curve in $\widehat{T/G}$ with the property $(P)$ and  $d_C \geq 3$. Then the naturally induced 1-form map from $H^0(T, \Omega_T )$ to $H^0(\widehat{W_C}, \Omega_{\widehat{W_C}})$ is an injective $G$-equivariant linear transformation. 
\end{lemma}
{\it Proof.} By choosing a fixed branched element of $\tau_C$ in (\req(tauC)) and with the translation of $T$ by a suitable  element in $T_G$, one may assume the dual of 1-form map, $H^0(\widehat{W_C}, \Omega_{\widehat{W_C}})^* \longrightarrow H^0(T, \Omega_T )^*$,  sends $H_1(\widehat{W_C}, \ZZ )$ to $H_1(T, \ZZ )$; hence one obtains the $G$-equivariant tori-homomorphism from the Jacobian of $\widehat{W_C}$ to $T$,
\be
\iota : J ( \widehat{W_C} ) \longrightarrow T \ .
\ele(iota)
Note that the dimension of $J ( \widehat{W_C} )$ is equal to $3d_C - 6 \geq 3$. It suffices to show the surjectivity of the above map $\iota$. Otherwise, by the $G$-eigenvalues of $H^0(T, \Omega_T )^*$, either the image or the cokernel of $\iota$ is a 1-torus with an order 7 automorphism group isomorphic to $G$, which is impossible as there is no 1-torus with an order 7 automorphism. Hence the result follows.    
$\Box$ \par \vspace{.2in} \noindent
For later use, we define the following biregular morphism group of $E(\omega)$,
\be
Aut(E(\omega), E(\omega)_{\langle m_\omega \rangle}) = \{ f : E(\omega) \stackrel{ \sim }{\longrightarrow} E(\omega)  | \ f (E(\omega)_{\langle m_\omega \rangle} ) = E(\omega)_{\langle m_\omega \rangle} \} ~ \simeq Aut(E(\omega)) \times \ZZ_3 \ ,
\ele(AEE)
where the above $\ZZ_3$ corresponds the translations by elements in $E(\omega)_{\langle m_\omega \rangle}$.

We now determine rational curves $C$ in $\widehat{T/G}$ with $d_C=3$,.
\begin{theorem} \label{thm:dC=3} For rational curves $C$ in $\widehat{T/G}$ with the property $(P)$ and $d_C=3$, 
the structure of $W_C$ and the number of $C$s, denoted by $\ell$, are given as follows.

(i) For $(T, G)= ( E(\omega)^3 , \langle m_\omega^3 \rangle )$ in  Example $1$, $W_C$ is the non-singular elliptic curve isomorphic to $E(\omega)$, and  $\ell = 513$.

(ii) For $(T, G)= ( A( \QZ (\mu )) , \langle m_\mu \rangle )$ in  Example $2$ , $W_C$ is the
genus 3 non-singular curve isomorphic to the Klein quartic $K$, and  $\ell = 14$. 
\end{theorem}
{\it Proof}. In this proof, $C$ always denotes a rational curve satisfying the assumption of this theorem. 

$(i)$. For a curve $C$, there exists at least one non-constant morphism $\rho_j$ among the three morphisms in (\req(rhoj)). Denote by $\ell^{(k)}, 1 \leq k \leq 3$, the number of $C$'s with exactly $k$ non-constant $\rho_j$s from $\widehat{W_C}$ onto $E(\omega)$.  By Lemma \ref{lem:Ed345}, $W_C$ is biregular to $E(\omega)$ via one of these $\rho_j$s. It is easy to see $\ell^{(1)} = 27$. For the computation of $\ell^{(2)}$, we consider a curve $C$ with exactly two non-constant $\rho_j$'s, say $\rho_1$ and $\rho_2$. Then by Lemma \ref{lem:Ed345}, there exists some $f \in Aut(E(\omega), E(\omega)_{\langle m_\omega \rangle})$ (defined in (\req(AEE))) such that $W_C$ is the graph $\Gamma_f$ in $E(\omega)^2$ (= the product of first two factors of $E(\omega)^3$). By $|Aut(E(\omega), E(\omega)_{\langle m_\omega \rangle})|= 18$, there are exactly 54 elliptic curves $W_C$  in $E(\omega)^3$ with $\rho_j$ constant only for $j=3$. Hence $\ell^{(2)}= 162$. By a similar argument, one has $\ell^{(3)}= 324$. Therefore $\ell = \sum_{k=1}^3 \ell^{(k)} = 513$.

$(ii)$. For $(T, G)= ( A( \QZ (\mu )) , \langle m_\mu \rangle )$, $\widehat{W_C}$ is a genus 3 Riemann surface. The morphism $\tau_C$ in  (\req(tauC)) defines a $G$-cover of $\widehat{W_C}$ over $C$ with three branched elements contained in $T_G$. We may assume $C = \PZ^1$ so that the branched loci of $\tau_C$ consists of 0, 1, and $\infty$. For an element $q \in T_G$, the $q$-translation $t_q$ of $T$, $ x \mapsto x + q $, commutes with $G$, hence $t_q(W_C)$ gives rise to another rational curve satisfying the assumption of this theorem. We  first consider the case when the identity element $o$ of $T$ corresponds a branched point of $\tau_C$. By Lemma \ref{lem:Tmu}, one has the $G$-equivariant isomorphism between $H^0(T, \Omega_T )$ to $H^0(\widehat{W_C}, \Omega_{\widehat{W_C}})$. Furthermore, by the bijective property of the morphism from $\widehat{W_C}$ onto $W_C$, the induced $G$-equivariant tori-morphism $\iota$ in (\req(iota)) is indeed an isomorphism compactible with the natural embedding of $\widehat{W_C}$ in $J(\widehat{W_C})$.  Therefore via $\iota$, $\widehat{W_C}$ is isomorphic to $W_C$, therefore $W_C$ is non-singular. By Lemma \ref{lem:g3G7}, $\widehat{W_C}$ is birational to the plane curve, $
y^7=x^2(x-1)$ ( hence birational to the Klein quartic $K$),
with $\tau_C$ corresponding to the $x$-projection. With a fixed identification of $J(K)$ with $T= A( \QZ (\mu ))$, the isomorphism $\iota$ belongs to the automorphism group generated by $G$ and the involution $i_T$ of $T$. Note that $G$ leaves $K$ invariant; while $i_T (K) \neq K$ since otherwise, the simple group $Aut(K)$ contains a non-trivial order 2 center element. Furthermore, by the same reason, $i_T (t_q (K)) \neq K$ for all $q \in T_G$. With the three choices for the branched element corresponding to $o$ of $T$, one obtains three embeddings of $K$ into $T$ passing through $o$. Indeed, we fix an arbitrary embedding $K$ into $T$ as $W_C$ with three branched points of $\tau_C$, denoted by $q_0(=o), q_1, q_2$. Then $q_j  \in T_G$ and the three embeddings of $K$ into $T$ are given by $t_{-q_j}(K)$ for $j=0,1,2$. As described in Example 2, $T_G$ is an order 7 (additive) subgroup of $T$. By which, one can show that $t_{-q_j}(K) \cap T_G $ are distinct subsets of $T_G$, hence $t_{-q_j}(K)$ three distinct curves in $T$ for $j=0,1,2$.
Therefore $t_{-q_j}(K), i_T(t_{-q_j}(K))$ for $j=0,1,2,$ are all $W_C (= \widehat{W_C})$ with $o$ as a branched point of $\tau_C$. By the translations $t_q$ for $q \in T_G$, one obtains all the $W_C$ in $T$, and the number $\ell = \frac{6 |T_G|}{3} = 14$. Then the result $(ii)$ follows.    
$\Box$ \par \vspace{.2in} \noindent

\section{Rational Curves $C$ in $\widehat{T/G}$ for $T= E(\omega)^3$ and $d_C \geq 4$ }
As in the previous section, $C$ denotes a rational curve in the rigid CY 3-fold $\widehat{T/G}$ with the property $(P)$. For $d_C \geq 4$, a natural problem would be the existence of $C$, the solution of which seems to be a difficult one, especially in the case $T=A( \QZ (\mu ))$.  
For $(T, G)$ in Example 1, no such $C$ exists for $d_C=4, 5$ by Lemma \ref{lem:Ed345}. Here we are going to investigate the cases for $d_C=6, 7$.

In this section, we shall denote $T= E(\omega)^3, G = \langle m_\omega^3 \rangle $. The elliptic curve $E(\omega)$ will always be represented as the non-singular compactification of the affine curve,
\be
y^3 = x (x-1) \ , \ \ (x, y ) \in \CZ^2 \ ,
\ele(cub)
with $dx/y^2$ as a base of holomorphic differentials, and the automorphism $m_\omega$ of $E(\omega)$ represented by the map, $(x, y) \mapsto (x , \omega y)$. The cover map of $E(\omega)$ over $E(\omega)/\ZZ_3 = \PZ^1$ corresponds the $x$-projection with the branched loci at $x=0, 1, \infty$, and the corresponding branched points in $E(\omega)$ will be denoted by $e_0, e_1, e_\infty$ respectively. The torus-involution of $E(\omega)$ corresponds to the map, $(x, y) \mapsto (\frac{1}{x}, \frac{-y}{x})$, with $e_1$ as the identity element of the order 3 (additive) subgroup $E(\omega)_{\langle m_\omega \rangle}(= \{ e_0, e_1, e_\infty \})$ of $E(\omega)$.   The $j$-th component of $E(\omega)^3$ will be denoted by $(x_j, y_j)$ satisfying the relation (\req(cub)) for $j=1,2, 3$. For a given $C$, one has the following commutative diagram,
$$
\begin{array}{lcl}
\widehat{W_C} & \stackrel{(\rho_1, \rho_2, \rho_3)}{\longrightarrow} & E(\omega)^3 \\
 \downarrow \tau_C & & \Pi \downarrow  \\
C & \stackrel{(\overline{\rho}_1, \overline{\rho}_2, \overline{\rho}_3)}{\longrightarrow} & ( \PZ^1)^3
\end{array}
$$
where $\tau_C$, $\rho_j$ are defined in (\req(tauC)), (\req(rhoj)) respectively, and  $\overline{\rho_j}$ is the $\ZZ_3$-quotient of $\rho_j$, and $\Pi$ is the triple-product of the covering map of $E(\omega)$ over $E(\omega)/\langle m_\omega \rangle = \PZ^1$. The morphism $\Pi$  sends $((x_j, y_j))_{j=1}^3$ to $(x_j)_{j=1}^3$, and it induces an one-one correspondence between elements of $(E(\omega)_{\langle m_\omega \rangle})^3$ and its $\Pi$-image set $(\{0, 1, \infty \})^3$. Note that $W_C$ could be singular with ${\rm Sing}(W_C) \subset (W_C)_G (= W_C \cap (E(\omega)_{\langle m_\omega \rangle})^3)$, and $(\rho_1, \rho_2, \rho_3)$ defines a bijective morphism from the normalization $\widehat{W_C}$ onto $W_C$. Hence $(\widehat{W_C})_G$, $\tau_C((\widehat{W_C})_G)$, $(W_C)_G$ and $\Pi((W_C)_G)$ are all bijective under the morphisms in the above diagram. 
\begin{lemma} \label{lem:6sm}
When $d_C=6$, the curve $W_C$ in $T$ is non-singular, (hence $\widehat{W_C} = W_C$).
\end{lemma}
{\it Proof}. Otherwise, there is an element $s \in (\widehat{W_C})_G$  corresponding to a singular point of $W_C$ in $T$.  Then the three positive integer, $r_{s,j}+1 \ ~ (j=1,2,3)$, in Lemma \ref{lem:lf3} must have at least one, say $r_{s,1}+1$, greater than or equal to 5 by (\req(gcdr)). The implies the degree $\delta$  of $\rho_1$ is greater than or equal to $5$ by setting $\rho = \rho_1$ in Lemma \ref{lem:rhoI}. By (\req(grho)), one has
$d_C = \delta + 2 +  \gamma^{(1)} \geq 7$, a contradiction to the assumption $d_C=6$. 
$\Box$ \par \vspace{.2in} \noindent
\begin{proposition} \label{prop:rhoS} (i) For $d_C \geq 6$, all $\rho_j$ are non-constant morphisms with $\delta_j := {\rm deg} (\rho_j) (= {\rm deg} (\overline{\rho}_j)) \geq 2$. The morphism $(\overline{\rho}_1, \overline{\rho}_2, \overline{\rho}_3)$ defines a bijective morphism between $C$ and its image $C^\dagger$ in $(\PZ^1)^3$, and elements of ${\rm Sing} (C^\dagger)$ and ${\rm Sing} (W_C)$ are in one-one correspondence under $\Pi$ with the isomorphic singularity-structure in their corresponding ambient $3$-spaces. Furthermore for $s \in (\widehat{W_C})_G$, the $\rho_j$-ramification index at $s$ is equal to the $\overline{\rho}_j$-ramification index at $\tau_C (s)$ for all $j$.  

(ii) Conversely, let ${\cal W}$ be a smooth curve with a $G$-action so that ${\cal W}/G \simeq \PZ^1$ and $|{\cal W}_G| = d \geq 6$, and $\vec{\rho} =(\rho_1, \rho_2, \rho_3)$ be a $G$-equivariant morphism of ${\cal W}$ into $E(\omega)^3$ under which ${\cal W}_G $ and $\vec{\rho}({\cal W}_G)$ are in one-one correspondence. Suppose that for each $s \in {\cal W}_G$, the local description of $\vec{\rho}$ near $s$ satisfies conditions in Lemma $\ref{lem:lf3}$, and the induced morphism 
$$
\vec{\overline{\rho}} =(\overline{\rho}_1, \overline{\rho}_2, \overline{\rho}_3): {\cal W}/G \longrightarrow (E(\omega)/\langle m_\omega \rangle )^3 = (\PZ^1)^3 
$$ 
defines a bijective morphism between ${\cal W}/G$ and $\vec{\overline{\rho}}({\cal W}/G)$, which is biregular outside $({\cal W}_G)/G$ with ${\rm deg} (\overline{\rho}_j) \geq 2$ for all $j$.  Then ${\cal W} = \widehat{W_C}$ for some $C$ in $\widehat{T/G}$ with the property $(P)$ and $d_C = d$. 
\end{proposition}
{\it Proof}. $(i)$. It is easy to see that a non-constant $\rho_j$ must have the degree $\geq 2$ ( since the degree-one case only occurs when $d_C=3$ by (\req(grho))), and there are at least two non-constant $\rho_j$'s, say $\rho_1, \rho_2$. 
One can also see that for $s \in (\widehat{W_C})_G$, the ramification indices of a non-constant $\rho_j$ at $s$ and $\overline{\rho}_j$ at $\tau_C (s)$ are the same. Indeed, if $(z_1, z_2 , z_3) = ( \alpha_1 t^{r_{s,1}+1} f_1(t^3) , ~ \alpha_2 t^{r_{s,2}+1} f_2(t^3), ~\alpha_3 t^{r_{s,3}+1} f_3(t^3) )$ is the local coordinate $t$-expression of $(\rho_1, \rho_2, \rho_3)$ centered at $s$ in Lemma \ref{lem:lf3}, then $\xi := t^3$ gives rise to the local coordinate in $C$ centered at $\tau_C (s)$, with the local description of $(\overline{\rho}_1, \overline{\rho}_2, \overline{\rho}_3)$ near $\tau_C (s)$ given by 
$$
\xi \mapsto ( \alpha_1^3 \xi^{r_{s,1}+1} f_1(\xi)^3 , ~ \alpha_2^3 \xi^{r_{s,2}+1} f_2(\xi)^3, ~\alpha_3^3 \xi^{r_{s,3}+1} f_3(\xi)^3 ) \ .
$$  
By simultaneously considering all non-constant $\rho_j$s as the $\rho$ in Lemma \ref{lem:rhoI}, one concludes that the morphism $\Pi$ sending $W_C$ to $C^\dagger$ is a local isomorphism at every element of $W_C$. Hence for the bijectivity property between $C$ and $C^\dagger$, one needs only to show that $(\overline{\rho}_1, \overline{\rho}_2, \overline{\rho}_3)$ defines an injective map from $C$ into $(\PZ^1)^3$. One may assume $C = \PZ^1$ with $B =  \{ a_j \in \CZ | 1 \leq j \leq d_C-1 \} \cup \{ \infty \} $ as the branched    
loci of $\tau_C$. Then $\widehat{W_C}$ can be represented by
$$
Y^3 = \prod_{j=1}^{d_C-1} (X - a_j )^{m_j} \ , \ \ (X, Y) \in \CZ^2 \ ,
$$ 
for some positive integers $m_j$s. The map $\tau_C$ corresponds to the $X$-projection, and $G$ is generated by $(X, Y) \mapsto (X, \omega Y)$. By Lemma \ref{lem:lf3}, the coordinate expressions for a non-constant $\rho_j$ and $\overline{\rho}_j$ have the following forms,
$$
\overline{\rho}_j : X \mapsto x_j = R_j(X) \ , \ \ \ ~ \rho_j: (X, Y) \mapsto (x_j , y_j)= ( R_j(X), Q_j(X)Y ) \ ,
$$  
where $R_j(X)$ is a degree-$\delta_j$ rational function of $X$ with $R_j^{-1} (\{0, 1 , \infty \}) = B$, and $Q_j(X)$ is rational function with $Q_j^{-1}(\{0, \infty \}) \subset B$. As the branched points of $\tau_C$, characterized by $Y= 0, \infty$, are one-one correspondence with $(W_C)_G$, the injectivity of $(\rho_1, \rho_2, \rho_3)$ is determined only by the $x_j$-values, hence equivalent to the injectivity of $(\overline{\rho}_1, \overline{\rho}_2, \overline{\rho}_3)$. This shows the bijectivity between $C$ and $C^\dagger$. 

It remains to show that all $\rho_j$ are non-constant. Otherwise, one has $ 6 \leq d_C \leq 9$, and we may assume that $(\overline{\rho}_1, \overline{\rho}_2)$ defines a bijective morphism between $C$ and $C^\dagger \subset (\PZ^1)^2$ with the isomorphic structure between ${\rm Sing} (C^\dagger) \subset (\PZ^1)^2$ and ${\rm Sing} (W_C) \subset E(\omega)^2$.  If $W_C$ is non-singular, one has $\PZ^1 = C \simeq C^\dagger$, which is a divisor in $(\PZ^1)^2$ linearly equivalent to $\delta_2 f_1 + \delta_1 f_2 $, where $f_j$ is a $j$-th fiber of $(\PZ^1)^2$. The adjunction formula for $C^\dagger$ implies
$$
\begin{array}{ll}
-2 & = c_1( ( \Omega_{(\PZ^1)^2} + [C^\dagger])_{C^\dagger}) = ( (\delta_2-2)f_1 + (\delta_1-2)f_2)\cdot (\delta_2 f_1 + \delta_1 f_2) \\
&= (\delta_2-2) \delta_1 + (\delta_1-2)\delta_2 \geq 0 \ , \ \ ~ ( {\rm by } \ \delta_j \geq 2) \ ,
\end{array}   
$$  
hence a contradiction.  Therefore $C^\dagger$ is singular, and $d_C \geq 7$ by Lemma \ref{lem:6sm}. By the results we have previously obtained, for $s \in B$ (= the branched loci  of $\tau_C$), the $\overline{\rho}_j$-ramification index of $s$  is $r_{s, j}+1$ for $j=1,2$. By (\req(grho)), one has $\delta_j = d_C -2 - \gamma^{(1)}_j$, hence $2 \leq \delta_j \leq 7$.  By Lemma \ref{lem:lf3}, $0 \leq r_{s, j} \leq 6$ and $r_{s, j} \neq 2, 5$ for  all $s, j$. For an element $s' \in B$ corresponding to a singular point of $C^\dagger$, one has $(r_{s', 1}, r_{s', 2}) = (1, 4), (4,1), (3, 6), (6,3)$. This implies $u (:= {\rm max}. \{r_{s, j} \ | \ s \in B, j=1, 2 \})$ is equal to $4$ or $6$. By Lemma \ref{lem:rhoI} and the $\ZZ_3$-quotient relation between $\rho_j$ and $\overline{\rho}_j$, one has $\overline{\rho}_j^{-1}(\{0, 1, \infty \}) = B$. Furthermore, by the property $C^\dagger \subset (\PZ^1)^2$, one has 
\be
|\overline{\rho}_j^{-1}(c)| \leq 3 \ , \ ~ {\rm for \ all } \ c \in \{0, 1, \infty \},  j=1,2 ,
\ele(2res)
 since otherwise, say  $|\overline{\rho}_1^{-1}(0)| \geq 4$, there are two elements in $\overline{\rho}_1^{-1}(0)$ with the same $\overline{\rho}_2$-value, a contradiction to the bijective relation between $C$ and $C^\dagger$. Let $s_0$ be an element in $B$ with $r_{s_0, j} = u$ for some $j$. Without loss of generality, we may assume $r_{s_0, 1} = u$ and $\overline{\rho}_1(s_0) = 0$. If $u=6$, one has $d_C=9$, $\delta_1 = 7$ and $\gamma^{(1)}_1 = 0$. By Lemma \ref{lem:rhoI}, we have $| \overline{\rho}_1^{-1}(0)|=1$. Hence $| \overline{\rho}_1^{-1}(1)| + | \overline{\rho}_1^{-1}(\infty)|=8$, which implies that either $| \overline{\rho}_1^{-1}(1)|$ or $| \overline{\rho}_1^{-1}(\infty)|$ is greater than 3, a contradiction to (\req(2res)). Therefore, $u $ must equal to $4$, hence $5 \leq \delta_1 \leq d_C-2 \leq 7$. By (\req(2res)) and Lemma \ref{lem:rhoI}, we have $| \overline{\rho}_1^{-1}(1)| = | \overline{\rho}_1^{-1}(\infty)|= 3$, and there are at least two elements $s \in \overline{\rho}_1^{-1}(c)$ for $c =1, \infty$ with $r_{s, 1}=1$ , which implies $r_{s, 2}=4$. Therefore there exists at least 4 elements $s \in B$ with $r_{s, 2}=4$, which is impossible by $\delta_2 \leq 7$. This completes the proof of $(i)$.  By reversing the arguments of $(i)$, one can easily see the converse statement $(ii)$ holds. 
$\Box$ \par \vspace{.2in} \noindent
{\bf Remark}. Note that the requirement in Proposition \ref{prop:rhoS} $(ii)$ on the local description in Lemma \ref{lem:lf3} near an element $s \in {\cal W}_G$ in ${\cal W}$  can be replaced by an equivalent one in ${\cal W}/G$. Hence in principal, the study of rational functions $\overline{\rho}_j$'s of $\PZ^1$ should provide an effective method of constructing rational curves $C$ in $\widehat{T/G}$ with $d_C \geq 6$. However, to obtain an explicit form of $\overline{\rho}_j$'s it remains a difficult task for the problem.  
$\Box$ \par \vspace{.2in} \noindent
For $d_C \geq 6$, by Proposition \ref{prop:rhoS} $(i)$, all $\rho_j$'s are non-constant morphisms. We shall denote the $\delta, \gamma^{(i)}$ in (\req(grho)) for $\rho= \rho_j$ by $\delta_j, \gamma^{(i)}_j$ for $j=1,2,3$ hereafter. Furthermore, by Lemma \ref{lem:rhoI} one has 
\be
\delta_j = \sum_{s \in \overline{\rho}^{-1}(z)} (r_{s, j}+1) \ \ , \ \ {\rm for} \ z = 0 , 1 , \infty \ , \ {\rm and} \ j = 1, 2, 3 \ .
\ele(deltaj)
(Note that we use the representation (\req(cub)) of $E(\omega)$).

 Now we consider the case $d_C=6$. The solutions of (\req(grho)) are: 
$$
(\gamma^{(0)}_j , \gamma^{(1)}_j ) = ( 0 , 2) \ , \ (3, 1) \ , \ (6 , 0 ) \ ; \ \ ~ \ ~  \delta_j = 4 -  \gamma^{(1)}_j \ .
$$
By Lemmas \ref{lem:lf3}, \ref{lem:rhoI} and \ref{lem:6sm}, $\gamma^{(1)}_j \neq 1$ and the ramification index $r_{s, j}+1$ for $s \in (W_C)_G$ is equal to $1$ or $4$. Hence $\delta_j = 2, 4$; furthermore one obtains the ramification indices for elements in $(W_C)_G$  by Lemma \ref{lem:lf3} as follows. When $\delta_j = 2$ ,  $|\rho_j^{-1}(e)| =2$ for all $e \in E(\omega)_{\ZZ_3}$, and $r_{s,j} = 0$   for all $s \in (W_C)_G $. When $\delta_j = 4$, there are two elements $e, e' \in E(\omega)_{\ZZ_3}$ with $|\rho_j^{-1}(e)|= |\rho_j^{-1}(e')|=1$, and the third one $e''$ with $|\rho_j^{-1}(e'')|=4$, hence $r_{s,j}=0$ for $s \in \rho_j^{-1}(e'')$ and the rest $r_{s,j} = 3$. The following two theorems give the full description on the structure of $W_C$ and its embedding in $E(\omega)^3$ for $d_C=6$.  
\begin{theorem} \label{thm:1d=6}
Assume $(T, G)= ( E(\omega)^3, \langle m_\omega^3 \rangle )$. For $d_C = 6$, we have $\delta_j = 2, 4$ for all $j=1,2,3$. If there is at least one $\delta_j=4$, the smooth curve $W_C$ is represented by  
\be
 Y^3 = X^4 (X-1)(X+1)(X-{\rm i})(X+{\rm i}) \ , \ \ (X, Y) \in \CZ^2 \ ,
\ele(64)
(i.e., $W_C$ is the compact smooth model of the above affine curve.) The embedding of $W_C$ in $E(\omega)^3$, unique up to permutations of components of $E(\omega)^3$, the symmetry of $Aut(W_C)$ on $W_C$ and  $Aut(E(\omega), E(\omega)_{\langle m_\omega \rangle})$ (defined in $(\req(AEE))$) on each factor $E(\omega)$ of $E(\omega)^3$, is given by the following standard representation.

(i). When $\delta_j =4$ for all $j$, the embedding of $W_C$ into $E(\omega)^3$ (with the coordinate $(\req(cub))$ of $E(\omega)$) is represented by
\bea(lll)
(x_1, y_1) = (X^4, Y), & (x_2, y_2)= (\frac{(X-1)^4}{(X+1)^4} , \frac{-2(X-1)Y}{X(X+1)^3} ) , &  (x_3, y_3) =  ( \frac{(X-{\rm i})^4}{(X+{\rm i})^4} , \frac{2{\rm i}(X-{\rm i})Y}{X(X+{\rm i})^3} ) . 
\elea(6444)

(ii). When exactly two $\delta_j$s equal to $4$, the embedding of $W_C$ into $E(\omega)^3$ is represented by
\bea(lll)
(x_1, y_1) = (X^4, Y), & (x_2, y_2)= (\frac{(X-1)^4}{(X+1)^4} , \frac{-2(X-1)Y}{X(X+1)^3} ) , & (x_3, y_3) =  ( \frac{ X(X+{\rm i})}{ {\rm i}(X-{\rm i})} , \frac{-Y}{X(X-{\rm i})} ) . 
\elea(6442)

(iii). When only one $\delta_j$ is equal to 4, the embedding of $W_C$ into $E(\omega)^3$ is represented by
\bea(lll)
(x_1, y_1) = (X^4, Y), & (x_2, y_2)= (\frac{X(X+1)}{X-1} , \frac{Y}{X(X-1)} ) , & (x_3, y_3) =  ( \frac{ X(X+{\rm i})}{ {\rm i}(X-{\rm i})} , \frac{-Y}{X(X-{\rm i})} ) .
\elea(6422)  
\end{theorem}
{\it Proof}. By the $X$-expression of $x_i$'s in each of  (\req(6444)), (\req(6442)) and (\req(6422)), one has an injective morphism of $X$ into $(\PZ^1)^3$, biregular outside $\{ 0, \infty, \pm 1, \pm {\rm i} \}$. By Proposition \ref{prop:rhoS} $(ii)$, these $(x_i, y_i)$'s expressions define the embedding of $W_C$ into $E(\omega)^3$ with the required properties. So we need only to derive the expression in each case.  We may assume $\delta_1=4$, and  $|\rho_1^{-1}(e_0)|= |\rho_1^{-1}(e_\infty)|=1$, $|\rho_1^{-1}(e_1)|=4$ (by a suitable $Aut(E(\omega), E(\omega)_{\langle m_\omega \rangle})$-transformation). Then $W_C$ can be represented by (\req(64)) with $\overline{\rho_1} : X \mapsto x_1=X^4$ and $\rho_1 : (X, Y) \mapsto (x_1, y_1)=(X^4, Y)$.
 The map $\tau_C$ of $W_C$ over $W_C/G = C = \PZ^1$  in (\req(tauC)) becomes the $X$-projection in the $(X, Y)$-representation of $W_C$, with the branched loci consisting of $X= 0, \infty, \pm 1, \pm {\rm i}$, and the corresponding branched elements in $W_C$ denoted by $s_0, s_\infty, s_{\pm 1}, s_{\pm {\rm i}}$ respectively. Then we have 
$$
\rho_1^{-1} (e_0) = \{ s_0 \} \ , \ ~ \rho_1^{-1} (e_\infty)= \{ s_\infty \} \ , \ ~ \rho_1^{-1} (e_1)= \{  s_{\pm 1} , s_{\pm {\rm i}} \} \ .
$$
Note that every birational morphism of $X \in \PZ^1$ preserving $\{ 0, \infty, \pm 1, \pm {\rm i} \}$ naturally induces an automorphism of $W_C$ commuting with $\tau_C$.
 
Claim: For $i= 2, 3$, the $\rho_i$-ramification index at $s_0, s_\infty$ is not equal to 4. Otherwise, without loss of generality we may assume $r_{s_\infty, 2}+1 = 4$, and $\overline{\rho}_2( \infty ) = \infty$. Then $\overline{\rho}_2$ is represented by $X \mapsto p(X) = \alpha (X- b )^4 $ for some $\alpha \in \CZ \setminus \{ 0 \}$ and $b \in \{ \pm 1 , \pm {\rm i} \}$. Then $p(0) = p( b')$ for $b' \in \{ \pm 1 , \pm {\rm i} \}$ and $b' \neq b$, which is impossible as there is no solution for such $p(X)$. 
     
$(i)$. We have $\delta_i= 4$ for $i=2, 3$. There are exactly two elements in $\{ \pm 1 , \pm {\rm i} \}$ with the $\overline{\rho}_i$-ramification index = 4, and $\overline{\rho}_i (0) = \overline{\rho}_i (\infty)$. We may assume $\overline{\rho}_i (0) = \overline{\rho}_i (\infty)=1$, and $\overline{\rho}_i$ is expressed by $X \mapsto x_i= \alpha \frac{X-b}{X-b'}$ for two distinct elements $ b, b' \in \{ \pm 1 , \pm {\rm i} \}$ and some $\alpha \in \CZ \setminus \{ 0 \}$. The solutions of $\{b , b' \}$ are $\pm 1, \pm {\rm i}$, hence we may assume $(b, b')=(1, -1), ({\rm i}, -{\rm i})$ by using the biregular morphism of (\req(64)) induced from the map $X \mapsto -X$ for the rest cases.   
As the rational map, $X \mapsto -{\rm i}X$, interchanges $\{ \pm 1 \}$ and $\{ \pm {\rm i} \}$, while fixing $0 , \infty$, one can write $\overline{\rho_2}, \overline{\rho_3}$ in the following forms after some suitable $Aut(E(\omega), E(\omega)_{\langle m_\omega \rangle})$-transformations:
$$
\begin{array}{ll}
\overline{\rho_2}(X) = \frac{(X-1)^4}{(X+1)^4} \ , & \overline{\rho_3}(X)= \frac{(X-{\rm i})^4}{(X+{\rm i})^4} \ . 
\end{array}
$$ 
Then the expression (\req(6444)) for the embedding $W_C$ in $E(\omega)^3$ immediately follows. And one has 
$$
\begin{array}{l}
\rho_2^{-1} (e_0) = \{ s_1 \} \ , \ ~ \rho_2^{-1} (e_\infty)= \{ s_{-1} \} \ , \ ~ \rho_2^{-1} (e_1)= \{ s_0 , s_\infty , s_{\pm {\rm i}} \} \ , \\
\rho_3^{-1} (e_0) = \{ s_{\rm i} \} \ , \ ~ \rho_3^{-1} (e_\infty)= \{ s_{-{\rm i}} \} \ , \ ~ \rho_3^{-1} (e_1)= \{ s_0 , s_\infty , s_{\pm 1} \} \ .
\end{array}
$$ 

$(ii)$. We may assume $\delta_2=4$ and $\delta_3=2$. By the argument in $(i)$, we have $(x_2, y_2)= (\frac{(X-1)^4}{(X+1)^4} , \frac{-2(X-1)Y}{X(X+1)^3} )$ with $
\rho_2^{-1} (e_0) = \{ s_1 \} $, $ \rho_2^{-1} (e_\infty)= \{ s_{-1} \}$, $ \rho_2^{-1} (e_1)= \{ s_0 , s_\infty , s_{\pm {\rm i}} \}$. The rational function $\overline{\rho}_3$ is of degree 2, unramified at $\{ \pm 1, \pm {\rm i}, 0 , \infty \}$. As the $(\overline{\rho}_1, \overline{\rho}_2)$-values of $\pm {\rm i}$ are both $(1, 1)$, by Proposition \ref{prop:rhoS} $(i)$, $\overline{\rho}_3(\pm {\rm i})$ are distinct with values in $\{ 0, \infty \}$; hence one may assume $\overline{\rho}_3({\rm i})= \infty$, $\overline{\rho}_3(-{\rm i})= 0$. Through the birational functions, $\pm \frac{X+1}{X-1}, \frac{-1}{X}$, one may further assume $\overline{\rho}_3(\infty)= \infty$. Then $\overline{\rho}_3^{-1}(0) = \{-{\rm i}, 0 \}, \{-{\rm i}, 1 \}, \{-{\rm i}, -1 \}$.  As there are two elements in $\{ \pm 1, \pm {\rm i}, 0 , \infty \} \setminus \overline{\rho}_3^{-1}(0, \infty)$ with the $\overline{\rho}_3$-value $= 1$, one concludes $\overline{\rho}_3^{-1}(0) = \{-{\rm i}, 0 \}$ with the relation $x_3 = \frac{ X(X+{\rm i})}{ {\rm i}(X-{\rm i})}$. The expression (\req(6442)) then follows, and we have 
$$
\rho_3^{-1} (e_0) = \{ s_0, s_{-{\rm i}} \} \ , \ ~ \rho_2^{-1} (e_\infty)= \{ s_\infty, s_{\rm i} \} \ , \ ~ \rho_3^{-1} (e_1)= \{ s_{\pm 1} \} \ .
$$

$(iii)$. We have $\delta_2= \delta_3 =2$. Both $\overline{\rho}_2$ and $\overline{\rho}_3$
are degree 2 rational functions of $X$ unramified at $\{ \pm 1, \pm {\rm i}, 0 , \infty \}$.  
By a suitable change of variables, one may assume that $\overline{\rho}_i$ takes different values at $0$ and $\infty$ for $i=2, 3$. Otherwise, say $\overline{\rho}_2(0)= \overline{\rho}_2(\infty) = \infty$, and $\overline{\rho}_2(1) = 0$. Then using (\req(deltaj)), one can show $\overline{\rho}_2(-1) = 1$. By applying the rational function $\frac{X+{\rm i}}{X-{\rm i}}$, the $\overline{\rho}_2$-values at $\infty, 0$ are now converted to $\pm 1$. Afrer a change of variables of $x_i$, we can further assume $\overline{\rho}_i(0)=0$, $\overline{\rho}_i(\infty) = \infty$ for $i=2, 3$. By changing $X$ to $-X$ if necessary, we may further assume $\overline{\rho}_2^{-1}(\infty) = \{ \infty , 1 \}$ and $\overline{\rho}_3^{-1}(\infty) = \{ \infty , {\rm i} \}$. This implies $x_2 = \frac{X(X+1)}{X-1}$, $x_3 = \frac{ X(X+{\rm i})}{ {\rm i}(X-{\rm i})}$ by Lemma \ref{lem:rhoI}. The expression (\req(6422)) then follows, and we have 
$$
\begin{array}{l}
\rho_2^{-1} (e_0) = \{ s_0 , s_{-1} \} \ , \ ~ \rho_2^{-1} (e_\infty)= \{s_\infty, s_1 \} \ , \ ~ \rho_2^{-1} (e_1)= \{ s_{\pm {\rm i}} \} \ , \\
\rho_3^{-1} (e_0) = \{ s_0, s_{-{\rm i}} \} \ , \ ~ \rho_2^{-1} (e_\infty)= \{ s_\infty, s_{\rm i} \} \ , \ ~ \rho_3^{-1} (e_1)= \{ s_{\pm 1} \} \ .
\end{array}
$$ 
$\Box$ \par \vspace{.2in} \noindent 
When all $\rho_j$s are of degree 2, we have the following result.
\begin{theorem} \label{thm:1d=62}
Assume $(T, G)= ( E(\omega)^3, \langle m_\omega^3 \rangle )$. For $d_C = 6$ and $\delta_j =2$ for all $j$, the embedding of $W_C$ into $E(\omega)^3$, unique up to permutations of components of $E(\omega)^3$, and $Aut(E(\omega), E(\omega)_{\langle m_\omega \rangle})$-symmetry on each factor of $E(\omega)^3$, is given by the following standard representation.

(i) When there exists a $\rho_j$ such that for the other two $k \neq j$, fibers of $\rho_k$ and $\rho_j$ over $E(\omega)_G$ coincide on one element, $W_C$ is described by the curve
\be 
Y^3 = X (X-1) (X+1) (3X-1) (3X+1) \ , \ ~ \ (X, Y) \in \CZ^2 \ ,
\ele(6a2i)  
and the embedding of $W_C$ into $E(\omega)^3$ is represented by
\bea(ll)
(x_1, y_1) = (\frac{(X-1)(3X+1)}{-4X}, \frac{Y}{2\sqrt[3]{2}X} ),& (x_2, y_2)= (\frac{-(X+1)}{(X-1)(3X+1)} , \frac{Y}{(X-1)(3X+1)} ) , \\
 (x_3, y_3) =  ( \frac{X-1}{(X+1)(3X-1)} , \frac{-Y}{(X+1)(3X-1)} ) .&
\elea(222i)

(ii) When there exists exactly only one pair, $\rho_j \neq \rho_k$, such that their fibers over $E(\omega)_G$ coincide on one element, $W_C$ is described by the curve
\be 
Y^3 = X (X-1)(X+1)(X-{\rm i})(X+{\rm i}) \ , \ ~ \ (X, Y) \in \CZ^2 \ ,
\ele(62)  
and the embedding of $W_C$ into $E(\omega)^3$ is represented by the following form:
\bea(lll)
(x_1, y_1) = (\frac{(X^2-1)}{{\rm i}(X+{\rm i})}, \frac{-Y}{X+{\rm i}} ),& (x_2, y_2)= (\frac{(X-{\rm i})}{{\rm i}(X^2-1)} , \frac{Y}{{\rm i}(X^2-1)} ) , &
 (x_3, y_3) =  ( \frac{-(X+1)}{X(X-1)} , \frac{Y}{X(X-1) }) .
\elea(222ii)
 
(iii) When any two of $\rho_j$'s has no common fibers over $E(\omega)_G$, the curve $W_C$ is represented by $(\req(62))$, and the embedding of $W_C$ into $E(\omega)^3$ is given by
\bea(lll)
(x_1, y_1) = (\frac{(X-1)(X-{\rm i})}{-2(1+{\rm i})X}, \frac{{\rm i}Y}{2X} ) , & (x_2, y_2)= (\frac{X(X+1)}{X-1} , \frac{Y}{X-1} ) , & (x_3, y_3) =  ( \frac{ X(X+{\rm i})}{{\rm i}(X-{\rm i})} , \frac{-Y}{X-{\rm i}} ) .
\elea(222iii)
\end{theorem}
{\it Proof}. By the $X$-expression of $x_i$'s in  (\req(222i)), (\req(222i)) and (\req(222iii)), one has an injective morphism of $X$ into $(\PZ^1)^3$ with different branched points in $X$-plane for the maps $x_1, x_2$ in each case. Hence by Proposition \ref{prop:rhoS} $(ii)$, the  $(x_i, y_i)$'s expression defines the embedding of $W_C$ into $E(\omega)^3$ with the required properties. So we need only to derive their expressions. 
As $\rho_j$ over $E(\omega)_G$ has the isomorphic structure as $\overline{\rho}_j$ over $\{ 0 , \infty, 1 \}$, the conditions in $(i)$-$(iii)$ can be carried over to those of $\overline{\rho}_j$'s. For a suitable coordinate of $C = \PZ^1$, we may assume the branched loci of $\tau_C$ consists of  six distinct elements, $0, \infty, 1, a, b, c$.

$(i)$. With suitable permutations and change of variables of $x_i's$,  we may assume 
\bea(lll)
\overline{\rho}_1^{-1}(\infty) = \{ 0 , \infty \}, & \overline{\rho}_1^{-1}(0) = \{ 1 , a \} &
\overline{\rho}_1^{-1}(1) = \{ b, c \} ; \\
\overline{\rho}_2^{-1}(\infty) = \{1, a \} & \overline{\rho}_2^{-1}(0) = \{\infty, b \} &
\overline{\rho}_2^{-1}(1) = \{0, c \} ,
\elea(rh12)
while for $\overline{\rho}_3$ one has 
$$
\overline{\rho}_3^{-1}(\infty) = \{b, c \}, \ ~ (\overline{\rho}_3^{-1}(0), \overline{\rho}_3^{-1}(1) )  = (\{\infty, 1 \},  \{0, a \} ) \ ~ {\rm or} \ (\{\infty, a \},  \{0, 1 \} ) \ .
$$ 
The solutions are given by $(a,b, c)= (\frac{-1}{3}, -1, \frac{1}{3}), (-3, 3, -1)$. The two cases are equivalent by the change of variables, $X \mapsto \frac{1}{X}$. So we may assume 
$(a,b, c)= (\frac{-1}{3}, -1, \frac{1}{3})$, hence $W_C$ is represented by (\req(6a2i)), and   
$$
\begin{array}{lll}
x_1 = \frac{(X-1)(3X+1)}{-4X} , & x_2 = \frac{-(X+1)}{(X-1)(3X+1)} , & x_3 = \frac{X-1}{(X+1)(3X-1)} \ .
\end{array}
$$
Then the expression (\req(222i)) follows.

$(ii)$. As in $(i)$, we may assume $\overline{\rho}_1$, $\overline{\rho}_2$ with the property (\req(rh12)). By the assumption of $(ii)$, the $\overline{\rho}_3$-values for $0$ and $\infty$ are distinct, the same for $1$ and $a$; furthermore the collection of these four values is equal to $\{ 0, \infty, 1 \}$. One may also assume $\overline{\rho}_3(\infty)=0$, $\overline{\rho}_3(0)= 1$, then we have the following descriptions  for $\overline{\rho}_3$-fibers:
$$
\begin{array}{lll}
(ii)_1: ~ \overline{\rho}_3^{-1}(0) =\{ \infty, a \} , & \overline{\rho}_3^{-1}(1) =\{ 0, b \} , & \overline{\rho}_3^{-1}(0) =\{ c, 1 \} ; \\
(ii)_2: ~ \overline{\rho}_3^{-1}(0) =\{ \infty, c \} , & \overline{\rho}_3^{-1}(1) =\{ 0, a \} , & \overline{\rho}_3^{-1}(0) =\{ b, 1 \}  ; \\
(ii)_3: ~ \overline{\rho}_3^{-1}(0) =\{ \infty, c \} , & \overline{\rho}_3^{-1}(1) =\{ 0, a \} , & \overline{\rho}_3^{-1}(0) =\{ b, 1 \}  ; \\
(ii)_4: ~ \overline{\rho}_3^{-1}(0) =\{ \infty, c \} , & \overline{\rho}_3^{-1}(1) =\{ 0, 1 \} , & \overline{\rho}_3^{-1}(0) =\{ a, b \} . \\
\end{array}
$$ 
The solutions of $(a, b, c)$ in the above cases are the same, given by $(a, b, c) = ({\rm i}, 1+{\rm i}, \frac{1+{\rm i}}{2}), (-{\rm i}, 1-{\rm i}, \frac{1-{\rm i}}{2})$. By changing the variable $X$ to $ \frac{1}{X}, \frac{X}{\rm i}, \frac{\rm i}{X}$, one can convert  $(ii)_1$ to  $(ii)_2, (ii)_3, (ii)_4$ respectively. So all cases are equivalent to those in $(ii)_1$. By  applying the changes variable, 
$$
\begin{array}{lll}
X \mapsto  \frac{(1-{\rm i})(X+{\rm i})}{2} , &(0, \infty, 1, -1, {\rm i}, -{\rm i}) \mapsto  (\frac{1+{\rm i}}{2}, \infty, 1, {\rm i}, 1+{\rm i}, 0) &   {\rm when } ~ (a, b, c) = ({\rm i}, 1+{\rm i}, \frac{1+{\rm i}}{2}) , \\
X \mapsto  \frac{2}{(1-{\rm i})(X+{\rm i})}, &(0, \infty, 1, -1, {\rm i}, -{\rm i}) \mapsto (1-{\rm i}, 0, 1, -{\rm i}, \frac{1-{\rm i}}{2}, \infty) & {\rm when} ~ (a, b, c) = (-{\rm i}, 1-{\rm i}, \frac{1-{\rm i}}{2}).
\end{array}
$$
the branched loci of $\tau_C$ becomes $\{ 0, \infty, \pm 1 , \pm {\rm i} \}$, hence $W_C$ is represented by (\req(62)). The two $X$-expressions of $x_i$'s derived from (\req(rh12)) and $(ii)_1$ are equivalent via the change of $X$-variables, $X \mapsto \frac{{\rm i}(X+{\rm i})}{X-{\rm i}}$. Hence we have $
x_1 = \frac{(X^2-1)}{{\rm i}(X+{\rm i})}$, $x_2 = \frac{(X-{\rm i})}{{\rm i}(X^2-1)}$, $   x_3 = \frac{-(X+1)}{X(X-1)}$. The expression (\req(222ii)) then follows.

$(iii)$. We may assume the $\overline{\rho}_1$-fibers over $\{ 0, \infty, 1 \}$ has the description in (\req(rh12)). With  suitable coordinates of $X$ and $x_2$,  $\overline{\rho}_2$-fibers are given by
$$
\overline{\rho}_2^{-1}(0) =\{0, c \} , ~ \overline{\rho}_2^{-1}(1) =\{ a, b  \} , ~ \overline{\rho}_2^{-1}(0) =\{ \infty,  1 \} \ .
$$
Then $\overline{\rho}_3$-fibers have the following expressions:
$$
\begin{array}{lll}
(iii)_1: ~ \overline{\rho}_3^{-1}(0) =\{0, b  \} , & \overline{\rho}_3^{-1}(1) =\{ 1, c  \} , & \overline{\rho}_3^{-1}(0) =\{ \infty , a  \} , \\
(iii)_2: ~ \overline{\rho}_3^{-1}(0) =\{0, a  \} , & \overline{\rho}_3^{-1}(1) =\{1, b  \} , & \overline{\rho}_3^{-1}(0) =\{ \infty , c \} , \\
(iii)_3: ~ \overline{\rho}_3^{-1}(0) =\{0, 1  \} , & \overline{\rho}_3^{-1}(1) =\{a, c  \} , & \overline{\rho}_3^{-1}(0) =\{ \infty , b \} , \\
(iii)_4:  ~ \overline{\rho}_3^{-1}(0) =\{0, a  \} , & \overline{\rho}_3^{-1}(1) =\{1, c  \} , & \overline{\rho}_3^{-1}(0) =\{ \infty , b \} . \\
\end{array}
$$
The solutions of $(a, b, c)$ are given by
$$
(a, b, c) = \left\{ \begin{array}{lll}
({\rm i}, -{\rm i}, -1) ,&(-{\rm i}, {\rm i}, -1) & {\rm for} \ (iii)_1 ; \\
(2, 1+{\rm i}, 1-{\rm i}) , &(2, 1-{\rm i}, 1+{\rm i}) & {\rm for} \ (iii)_2 ; \\
({\rm i}, \frac{1+\rm i}{2}, 1+{\rm i}) , & (-{\rm i}, \frac{1-\rm i}{2}, 1-{\rm i}) & {\rm for} \ (iii)_3 ; \\
(\frac{b^2}{b^2-b+1}, b , \frac{b^2}{b^2-b+1}) & & {\rm for} \ (iii)_4 \ .
\end{array}\right.
$$
For the solutions of $(iii)_4$, the injective property of $(\overline{\rho}_1, \overline{\rho}_2, \overline{\rho}_3)$ in Proposition \ref{prop:rhoS} $(i)$ implies $b=1$, a contradiction to the choice of $ b$, hence $(iii)_4$ can be ruled out. For each one of the rest cases, the structures of $(\overline{\rho}_1, \overline{\rho}_2, \overline{\rho}_3)$ provided by the two solutions of $(a, b, c)$ are equivalent under the symmetries of this theorem after a suitable coordinate-change of $X$. Indeed for $(iii)_1$, by  $X \mapsto \frac{-(X-{\rm i})}{X+{\rm i}}$, one can transform the structure of $\overline{\rho}_j$'s of the second solution to the first one. The same is for $(iii)_2$ by $X \mapsto \frac{2}{X}$, and for $(iii)_3$ by $X \mapsto \frac{1}{X}$. Hence it suffices to consider the structures given by the first solution of $(iii)_i$ for $i=1, 2, 3$. By the maps, $X \mapsto 1-X$ or $1-(1-{\rm i})X$, one can reduce the first solution of $(iii)_2$, $(iii)_3$ to that of $(iii)_1$. Hence the  structure of $\overline{\rho}_j$'s is uniquely given by 
$$
\begin{array}{lll}
x_1= \frac{(X-1)(X-{\rm i})}{-2(1+{\rm i})X}, & x_2 = \frac{X(X+1)}{X-1}  , & x_3 = \frac{ X(X+{\rm i})}{{\rm i}(X-{\rm i})} ,
\end{array}
$$
with the branched loci of $\tau_C$ = $\{ 0, \infty, \pm 1 , \pm {\rm i} \}$. The results of $(iii)$ then immediately follow.
$\Box$ \par \vspace{.2in} \noindent
{\bf Remark}. Note that the involved symmetries in Theorem \ref{thm:1d=6} are more than those in Theorem \ref{thm:1d=62} ( with the extra $Aut(W_C)$-symmetries on $W_C$ for the former one). However, the $W_C$'s in (\req(64)), (\req(6a2i)), (\req(62)) pairwise are not $\ZZ_3$-equivariant. The counting of curves $C$ could be obtained, though cumbersome but workable, for each cases in Theorems \ref{thm:1d=6} and \ref{thm:1d=62}. 
$\Box$ \par \vspace{.2in} 
Now we consider the case $d_C=7$. The curve $W_C$ is singular by Lemma \ref{lem:Ed345} $(iii)$. Using the relation (\req(grho)), $(\delta_j, \gamma^{(0)}_j, \gamma^{(1)}_j)$ takes one of the values $(3, 2, 2), (4, 5, 1), (5, 8, 0)$ for all $j$. Hence the singular points of $W_C$  correspond to those element $s$ in $ (\widehat{W_C})_G$ with $r_{s, j}=4$ and $r_{s, k} =1$ for some $j \neq k$ by Lemma \ref{lem:lf3}. Furthermore, all $r_{s, j}+1$ for $s \in (\widehat{W_C})_G$ are positive integers $\leq 5$ and $\neq 3$.  Denote the branched loci of $\tau_C$ in $C$ by $S$, which can be identified with $(\widehat{W_C})_G$ via $\tau_C$, hence $|S| = 7$. Without loss of generality, we may assume $\delta_1=5$ and $r_{s, 1}=4$ for some $s \in S$. 

Claim: $|{\rm Sing} (C^\dagger)| (= |{\rm Sing} (W_C)| ) = 2$, and for suitable $\PZ^1$-coordinates $X$ for $C$, one has $S = \{ 0, \infty \} \cup \{ \alpha \in \CZ \ | \ \alpha^5 =1 \}$ with  
${\rm Sing} (C^\dagger)$ corresponding to $\{ 0 , \infty \}$ and the expression of $\overline{\rho}_1$ given by $x_1 = X^5$. In fact, the $X$-coordinates of $C$ can be made so that $\infty \in S$ with $r_{\infty, 1} = 4$. Then $\overline{\rho}_1^{-1}(\infty) =\{ \infty \}$. By (\req(deltaj)), the $\overline{\rho}_1$-fibers over $\{0, 1 \}$ are described by one of the following list for a suitable $x_1$-coordinate,
$$
\begin{array}{lll}
7_{(i)}: & | \overline{\rho}_1^{-1}(0) | = 1 , & | \overline{\rho}_1^{-1}(1) | = 5 ; \\
7_{(ii)}: & | \overline{\rho}_1^{-1}(0) | = 2 , & | \overline{\rho}_1^{-1}(1) | = 4 ; \\
7_{(iii)}: & | \overline{\rho}_1^{-1}(0) | = 3 , & | \overline{\rho}_1^{-1}(1) | = 3 , \
\end{array}
$$ 
with $\overline{\rho}_1$-ramification indices of elements in $S$ and the $X$-expression 
of $x_1$ for a suitable coordinate $X$ of $C$ as follows:
$$
\begin{array}{llll}
7_{(i)}: & x_1 = X^5 , & S = \{ 0, \infty \} \cup \{ \alpha  \ | \ \alpha^5 =1 \} , & {\rm where} \ r_{0, 1} = 4 \ , \ r_{\alpha , 1} = 0  ; \\
7_{(ii)}: &  x_1 = \frac{X^4 (X-1)}{a^4(a-1)} , & S = \{ 0, \infty , 1 , a, b_1, b_2, b_3 \} , & {\rm where} \ r_{a,1} = 1, \ r_{b_j, 1} = 0 \ , (j=1,2,3) \ ; \\
7_{(iii)}:& x_1 = \frac{X(X-a_1)(X-a_2)}{(1-a_1)(1-a_2)}  , & S = \{ 0, \infty , 1 , a_1, a_2, a_3, a_4 \} & {\rm where} \ r_{0,1} = r_{1,1} = 0 , \ r_{a_j,1} = 1, (1 \leq j \leq 4 ) .
\end{array}
$$
In the case $7_{(iii)}$, the singularities of $C^\dagger$ consists of 5 elements, $\infty, a_j (1 \leq j \leq 5)$, among which there are at two elements with $\overline{\rho}_k$-ramification index 5 for $k=2$ or $3$ by the property (\req(gcdr)). Then the structure of $\overline{\rho}_k$ is equivalent to the case $7_{(i)}$, which implies $|{\rm Sing } (C^\dagger)|=2$, hence a contradiction. In the case $7_{(ii)}$, $a$ and $b_j$'s satisfy the relation
$$
\frac{d x_1}{dX} (a ) = 0 \ , \ \ x_1(a)= x_1(b_j) \ , \ \ j=1, 2, 3 \ .
$$   
By computations, there are no solution of $a$ and $b_j$'s so that $S$ consists of 7 distinct elements. Therefore $7_{(i)}$ is the only description of $\overline{\rho}_1$, hence follows the claimed results. 

Now we make the identification $S = \{ 0, \infty \} \cup \{ \alpha \in \CZ \ | \ \alpha^5 =1 \}$ with ${\rm Sing} (C^\dagger)$ corresponding to $\{ 0 , \infty \}$. The degree $\delta_2$ is equal to $3$ or $4$, since otherwise, the map $\overline{\rho}_2$ is equivalent to the case $7_{(i)}$. Then the five elements, $ \alpha \in \CZ $ with $ \alpha^5 =1$, take the same value of $(x_1, x_2)$, hence can't be mapped injectively in $(\PZ^1)^3$. When $\delta_2 = 3$, by (\req(deltaj)), one may assume the $X$-expression of $x_2$ is given by $x_2 = (d+1) \frac{X^2(X-d)}{X-1}$ for some primitive $5$th root of unity $d$, satisfying the relation $x_2 (d^2) = x_2(d^3) = x_2(d^4) = 1$, which is impossible. Therefore $\delta_2 = 4$, then one may write the $X$-expression of $x_2$ either as $x_2 = \frac{\beta^2}{(\beta -1)^4} \frac{(X-1)^4}{X^2}$ with $x_2(\beta^k)= 1$ for $1 \leq k \leq 4$, or $x_2 = \beta^{\ell+1}
 \frac{(X-1)^4}{(X- \beta)(X-\beta^\ell)}$ for some $1 \leq \ell \leq 4$ with $\frac{d x_2}{dX} (0) = 0$ and $x_2(0)= x_2(\beta^k) = 1$ for $ 2 \leq k \neq \ell \leq 4$, where $\beta$ is a primitive $5$th root of unity. One can easily see that there is no such $X$-expression of $x_2$. Therefore, we have shown the following result: 
\begin{proposition} \label{prop:d=7}
For $(T, G)= ( E(\omega)^3, \langle m_\omega^3 \rangle )$, there is no rational curve $C$ in $\widehat{T/G}$ with $d_C = 7$.
\end{proposition}

\section{ Conclusions and Perspectives }
In this work, we have studied the structure of CY 3-fold $\widehat{T/G}$ associated to 3-torus-orbifolds $T/G$ with only isolated singularities, then conduct a primitive investigation of a  rational-curve problem on $\widehat{T/G}$. By the theory of CM-type abelian varieties and a topological relation on the global structure of $T$ with its $G$-fixed pont set $T_G$, the isomorphic classes of $(T, G)$ are determined, which are represented  by $( E(\omega )^3, \ZZ_3)$ in Example 1, or $(A(\QZ(\mu)), \ZZ_7)$ in Example 2  in Sect. 2. For each case, the crepant resolution $\widehat{T/G}$ is a rigid CY 3-fold with a simple structure of exceptional divisors. By studying some special $G$-invariant curves in the 3-torus $T$, we investigate the rational-curve problem about a $\PZ^1$-curve $C$ in  $\widehat{T/G}$ not in exceptional divisors with the property $(P)$ (in Sect. 3), then introduce the counting number $d_C$ of elements in $C$ intersecting the exceptional set of $\widehat{T/G}$. When $d_C$ equals to the minimal number 3, we have obtained  the structure of  the corresponding $G$-invariant curve $W_C$ in $T$, and the total number of such rational curve $C$ in $\widehat{T/G}$. For $d_C \geq 4$, progresses on this rational-curve problem in $\widehat{T/G}$ are made only in the case $(T, G)=(E(\omega )^3, \ZZ_3)$, where a method was developed by reducing the problem of $C$ in the CY space to one about rational curves in $(\PZ^1)^3$. We have shown that for a given $d_C$, the existence of $C$ appears to be a non-trivial problem by examining each case for $d_C \leq 7$. Only when $d_C=6$, one has the positive solution for the problem by performing explicit polynomial-calculations. Indeed, the complete characterization of the $G$-invariant curve $W_C$ in $T$ has been derived in the case $d_C=6$. The results obtained so far seems to suggest some reasons skeptical to the existence of rational curve $C$ with higher $d_C$, even though the justification has not been done yet. The qualitative understanding of the $G$-invariant curves $W_C$ in the 3-torus $T$ would still be required for further study on this rational curve problem in these rigid CY 3-folds $\widehat{T/G}$. Recently, progresses on the modularity property of certain rigid CY 3-folds defined over number fields have been made in arithmetic geometry (see, e.g. \cite{Y} and reference therein). The rigid CY 3-folds $\widehat{T/G}$ studied in this paper appears to be among interesting models deserving for close inspection to the possible modularity feature.
Owing to the intricacy of the defining equations of $T/G$, much remains to be discovered in the arithmetic aspect of $\widehat{T/G}$.

\section*{Appendix: 3-tori $T$ with a finite cyclic automorphism group $G$ and $|T_G | < \infty$ }
In this appendix, we are going to derive the complete list of $(T, G)$ for a 3-torus $T$ and a finite cyclic automorphism group $G \subset Aut(T)$ with $|T_G | < \infty$. 

First we state the following algebraic fact, whose proof can be found in \cite{CW}.
\begin{lemma} \label{lem:QG}
Let $M$ be a finite-dimensional vector space over $\QZ$, and $G$ a finite abelian subgroup of $Aut_{\QZ} ( M ) \subset End_{\QZ} ( M )$ with $|G| = d$. Denote $\QZ G$ the group ring of $G$ over $\QZ$, regarded as a subring of $End_{\QZ} ( M )$. Suppose $M$ is a simple $\QZ G$-module. Then $G$ is a cyclic group, and $\QZ G = \QZ ( \zeta_d )$, ${\rm dim}_{\QZ} M = \phi (d)$, where   $\phi$ is the Euler function, and $\QZ ( \zeta_d)$ is the cyclotomic field of degree $d$ with $\zeta_d = e^{\frac{2 \pi {\rm i}}{d}}$. 
\end{lemma}
$\Box$ \par \vspace{.2in}\noindent
\begin{proposition} \label{prop:CM}
Let $T$ be a $n$-torus ($n \geq 2$), $G$ an order $d$ cyclic subgroup of $Aut(T)$ with $|T_G| < \infty$, and $\theta$ a generator of $G$. Let $\{ \lambda_j \}_{j=1}^s$ be the collection of all distinct eigenvalues of the differential $(d \theta)_o$ at $o \in T$. Denote $\Phi = \{ \varphi_j \}_{j=1}^s $ where $\varphi_j$ is the field-embedding of $\QZ( \zeta_d )$ into $\CZ$ with $\varphi_j (\zeta_d) = \lambda_j $, and $h_d :=$ the ideal-class number of $\QZ(\zeta_d)$. Then 

(i) $2n$ is divisible by $\phi (d)$, and $s \geq \frac{\phi(d)}{2}$. 

(ii) When $2n = \phi (d)$, (hence $n= s = \frac{\phi(d)}{2}$ by (i)), the $n$-torus $T$ is an abelian variety of type $(\QZ(\zeta_d), \Phi)$ in the sense of $\cite{ST}$ under the identification of $m_{\zeta_d}$ (the $\zeta_d$-multiplication) with $\theta$ in $End_{\QZ} (T)$. Furthermore, for a given $\Phi$, the number of (non-isomorphic) abelian varieties of type $(\QZ(\zeta_d), \Phi)$ is equal to $h_d$. In particular, when $h_d = 1$, $T$ is isomorphic to $A(\QZ(\zeta_d), \Phi)$ (:= the abelian variety $(\QZ(\zeta_d)\otimes_{\QZ} \RZ )/ (\ZZ(\zeta_d)\otimes 1)$ with the natural induced complex structure).

(iii) When $s = \frac{\phi(d)}{2}$,  $( \QZ( \zeta_d) , \Phi )$ is a CM-field. Denote $k= \frac{2n}{\phi(d)}$. Then $T$ is isomorphic to $k$-product of the $( \QZ( \zeta_d) , \Phi )$-type abelian variety , and $\theta$ can be chosen to correspond to the $k$-product of $m_{\zeta_d}$. Furthermore, for a given $\Phi$, there are exactly ${k + h_d - 1 \choose k}$ such (non-isomorphic) $n$-dimensional abelian varieties $T$. In particular, $T$ is isomorphic to $( A(\QZ(\zeta_d), \Phi))^k$ when $h_d=1$.
\end{proposition}
{\it Proof}. $(i)$. Consider the induced $G$-action of the $2n$-dimensional $\QZ$-space $H_1 (T, \QZ)$.   By $|T_G| < \infty$, $G$ can be regarded as a subgroup of $Aut_{\QZ} ( H_1 (T, \QZ) )$. Hence there is a simple $\QZ G$-submodule $M$ of $H_1 (T, \QZ)$. By Lemma \ref{lem:QG}, we have $\QZ G = \QZ ( \zeta_d )$, so $H_1 (T, \QZ)$ becomes a $\QZ ( \zeta_d )$-vector space, which implies $ \phi(d) | 2n$. By the identification of $H_1 (T, \QZ) \otimes_{\QZ} \CZ$ with $H^{1,0}(T)^* \oplus H^{0, 1}(T)^*$, and $H^{1,0}(T)^*$ with the tangent space of $T$ at $o$, on which there are $s$-distinct eigenvalues for the linear map $(d \theta)_o$, one has the $\CZ$-dimension of 
$ \QZ ( \zeta_d )\otimes_{\QZ} \CZ$ is less than or equal to $2s$. Hence $s \geq \frac{\phi(d)}{2}$. 

$(ii)$. In the case $2n = \phi (d)$, one has the isomorphism of $\QZ (\zeta_d)$ into $End_{\QZ} (T)$, which sends $\zeta_d$ to $\theta$. Then $T$ is abelian variety of type $(\QZ(\zeta_d), \Phi)$ \cite{Sh}. The rest of conclusions in $(ii)$ follows from  Proposition 17 in \cite{ST}.

$(iii)$. By $s = \frac{\phi(d)}{2}$, any simple $\QZ G$-submoulde of $H_1 (T, \QZ)$ can be identified with $\QZ ( \zeta_d )$, then follows the CM-field $( \QZ( \zeta_d) , \Phi )$. Furthermore, with the $\theta$-induced action, $H_1 (T, \ZZ)$ becomes a $\ZZ [\zeta_d]$-module. 
By the classification of regular $\ZZ [\zeta_d]$-modules \cite{K, Re}, $H_1 (T, \ZZ)$ is isomorphic to a direct sum of ideals, $I_j (1 \leq j \leq k)$, in $\ZZ [\zeta_d]$,
$$
H_1 (T, \ZZ) \ \simeq I_1 + \cdots + I_k \  \ ~ {\rm as} \ \ZZ [\zeta_d]-{\rm modules} \ .
$$
Note that the rank $k$ and the ideal-classes of $I_j$ are the only invariants for the above $\ZZ [\zeta_d]$-module structure. This implies the number of the $\ZZ [\zeta_d]$-modules of rank $k$ is equal to ${k + h_d - 1 \choose k}$.  By "$\otimes_{\ZZ} \QZ$" (or "$\otimes_{\ZZ} \RZ$") with the above relation, one obtains the $\QZ(\zeta_d)$ ($\QZ(\zeta_d)\otimes_{\QZ} \RZ$ resp.) structure of $H_1 (T, \QZ)$ ( $H_1 (T, \RZ)$ resp.), which are compatible with the above $\ZZ [\zeta_d]$-module decomposition of $H_1 (T, \ZZ)$. Hence as real-tori, 
\be
T = H_1 (T, \RZ)/ H_1 (T, \ZZ) ~ \simeq \prod_{j=1}^k ( \QZ(\zeta_d)\otimes_{\QZ} \RZ)/(I_j \times 1) \ ,
\ele(TIj)
by which the $\theta$-action on $T$ corresponds to the product of $m_{\zeta_d}$ of the factors. In order to identify $T$ with the $k$-product of the $( \QZ( \zeta_d) , \Phi )$-type abelian varieties through (\req(TIj)), one needs to show that (\req(TIj)) gives rise to an isomorphism of complex-tori by using $\Phi$. For this purpose, it suffices to show the multiplications on $H_1(T, \RZ)$(= the tangent space of $V$ at $o$) by complex numbers are elements in the algebra over $\RZ$ generated by $(d \theta)_o$. By $s = \frac{\phi(d)}{2}$, for $c \in \CZ$, there exist $r_{s,j} \in \RZ$ ($0 \leq j \leq \phi (d)-1$) such that 
$$
\left(\begin{array}{c} c \\ \vdots \\ c \end{array} \right) = \sum_{j=0}^{\phi(d)-1} r_{s,j} \left(\begin{array}{c} \lambda_1^j \\ \vdots \\ \lambda_s^j \end{array} \right) \  , \ ~ ~  \ s = \frac{\phi(d)}{2} \ . 
$$
By which, through the field-embeddings in $\Phi$, one concludes that the $c$-multiplication is expressed by $\sum_{j=0}^{\phi(d)-1} r_{s,j} (d \theta)_o^j$. The results of $(iii)$ then follow.  
$\Box$ \par \vspace{.1in}\noindent
{\bf Remark}. The statement $(ii)$ is a special case of $(iii)$ for $k=1$.
And in the assumption of the above proposition, the "cyclic" condition of $G$ can be replaced by "abelian" one, but with the same conclusions, by Lemma \ref{lem:QG}.
$\Box$ \par \vspace{.2in}\noindent
For a cyclotomic field $\QZ(\zeta_d)$ of degree $d$, and a positive integer $j (<d)$ relatively prime to $d$, there is the field-embedding of $\QZ(\zeta_d)$ into $\CZ$, denoted by $\psi_j (= \psi_{d, j})$, defined by $\psi_j ( \zeta_d ) = \zeta_d^j$. 
\begin{proposition} \label{prop:TG3}
Let $T$ be a 3-torus, and $G$ a cyclic subgroup of $Aut(T)$ with $|T_G| < \infty$ and $d (=|G|) \geq 3$. A list of $(T, G)$ (up to isomorphisms) with the description of a generator of $G$ is given as follows:
$$
\begin{array}{| c |c | c| }
\hline 
|G| & T & \theta ~  ( {\rm a \ generator \ of } \ G ) \\
\hline 
3 & E(\omega)\times E(\omega)\times E(\omega)  & m_\omega \times m_\omega \times m_\omega  \\
\hline 
3 & * & (d \theta)_o \sim {\rm dia} ~ [\omega, \omega, \omega^2 ] \\
\hline 
4 & E({\rm i})\times E({\rm i})\times E({\rm i})  & m_{\rm i} \times m_{\rm i} \times m_{\rm i}  \\
\hline 
4 & ** & (d \theta)_o \sim {\rm dia} ~ [{\rm i}, {\rm i}, -{\rm i} ] \\
\hline 
6 & E(\omega)\times E(\omega)\times E(\omega)  & -m_\omega \times -m_\omega \times -m_\omega  \\
\hline 
6 & * & (d \theta)_o \sim {\rm dia} ~ [-\omega, -\omega, -\omega^2 ]\\
\hline 
7 & A(\QZ (\zeta_7) , \{\psi_1, \psi_2, \psi_3 \}) , ~ A(\QZ (\zeta_7) , \{\psi_1, \psi_2, \psi_4 \}) &  m_{{\zeta_7}} \\
\hline 
9 & A(\QZ (\zeta_9) , \{\psi_1, \psi_2, \psi_4 \}) , ~ A(\QZ (\zeta_9) , \{\psi_1, \psi_4, \psi_7 \}) &  m_{{\zeta_9}} \\
\hline 
14 & A(\QZ (\zeta_7) , \{\psi_1, \psi_2, \psi_3 \}) , ~ A(\QZ (\zeta_7) , \{\psi_1, \psi_2, \psi_4 \}) &  - m_{{\zeta_7}} \\
\hline 
18 & A(\QZ (\zeta_9) , \{\psi_1, \psi_2, \psi_4 \}) , ~ A(\QZ (\zeta_9) , \{\psi_1, \psi_4, \psi_7 \}) &  - m_{{\zeta_9}} \\
\hline 
\end{array}
$$
where $E(\omega), E({\rm i})$ are the 1-tori with $(E(\omega),  m_\omega)$ defined in {\rm Example 1}, and  $E({\rm i})= \CZ/(\ZZ + \ZZ {\rm i})$ , $m_{\rm i}: [z] \mapsto [{\rm i}z]$. The structure of classes $"*"$, $"**"$ has yet completely determined, and $"\sim"$ means the conjugation of matrices.

As a consequence, when $G \subset Saut (T)$, $(T, G)$ is isomorphic to either $(E(\omega)^3, \langle m_\omega^3 \rangle )$ (in {\rm Example 1}), or $((\QZ (\zeta_7) , \{\psi_1, \psi_2, \psi_4 \}) , \langle m_{{\zeta_7}} \rangle)$ (in {\rm Example 2}).
\end{proposition}
{\it Proof}. By Proposition \ref{prop:CM} $(i)$, we have $\phi(d) | 6$. Using the expression of $\phi(d)$ in terms of its prime factors $p$, $\phi (d) = d \prod_{p | d} (1 - \frac{1}{p})$, one obtains the conclusion: $d = 3, 4, 6 , 7, 9 , 14 , 18$, hence the 
ideal-class number $h_d$ of $\QZ( \zeta_d )$=1.  When $d = 7, 9 , 14 , 18$, we have $\phi(d)=6$. Then the classification in those cases follows from Proposition \ref{prop:CM} $(ii)$, and $\zeta_{14} = - \zeta_7$, $\zeta_{18} = - \zeta_9$. When $d= 3, 4, 6$, we have $\phi (d) = 2$. By $|T_G | < \infty$, any eigenvalue of $(d \theta)_o$ must be a primitive $d$-th root of unity. Hence $s$ (in Proposition \ref{prop:CM}) is equal to 1 or 2. When $s=1$, we obtain the structure of the 3-product of 1-torus in the above table by Proposition \ref{prop:CM} $(iii)$. When $s=2$,    the eigenvalues $(d \theta)_o$ can be determined as in the table, while the tori-structure of  $T$ remains to be determined.   
$\Box$ \par \vspace{.1in}\noindent
{\bf Remark}. In the table of the above proposition,  the detailed structures of classes $"*"$ and $"**"$ haven't been identified yet. However, these classes do not possess the property $G \subset Saut(T)$, of which we required in the main content of this paper. But there does exist some torus in $"*"$ and $"**"$. Indeed, one can easily see that $E(\omega)^3$ is an example in $"*"$. For the class $"**"$, it contains the Jacobian of the following one-parameter family of genus-3 hyperelliptic curves, appeared in the study of chiral Potts $N$-state model for $N=4$ (see, e.g., \cite{B91}):
$$
t^4 = \frac{(1 - \kappa
\lambda)(1 - \kappa \lambda^{-1})}{1 - \kappa^2 } \ , \ \ \ (t, \lambda )
\in \CZ^2 \ ,
$$
for $\kappa \in \CZ \setminus \{0, \pm 1 \}$. 
$\Box$ \par \vspace{.2in}\noindent

\section*{Acknowledgements}
The author is grateful to N. Yui for the friendly hospitality during the "BIRS Workshop: Calabi-Yau Varieties and Mirror Symmetry" December 6-11, 2003, in Banff Centre, Alberta, Canada. He would also like to thank K. Oguiso for the informative discussions on his work related to some topic of this note during the BIRS workshop.  
This work has been supported by NSC 92-2115-M-001-023.

\end{document}